\documentclass[11pt,a4paper]{amsart}

\usepackage{amsmath,amstext,amssymb,amsthm,amsfonts,xy}

\theoremstyle{plain}
\newtheorem*{theorem*}{Theorem}
\newtheorem*{remark*}{Remark}
\newtheorem*{example*}{Example}

\newtheorem*{conjecture*}{Conjecture}

\theoremstyle{definition}

\theoremstyle{remark}

\newcommand{\bR}{{\mathbb R}}

\def\quotient#1#2{%
    \raise1ex\hbox{$#1$}\Big/\lower1ex\hbox{$#2$}%
}

\input xy
\xyoption{all}

\begin{document}

\date{January 25, 2017}

\setcounter{tocdepth}{1}
\title[Schwartz functions on real algebraic varieties]
{Schwartz functions on real algebraic varieties}

\author[Boaz Elazar and Ary Shaviv]{Boaz Elazar$^\dag$ and Ary Shaviv$^\dag$}

\address[]{Dept. of Mathematics, The Weizmann Institute of Science,
Rehovot 76100, Israel}
\email{boaz.elazar@weizmann.ac.il, ary.shaviv@weizmann.ac.il}

\thanks{$^\dag$Supported in part by ERC StG grant 637912 and ISF grant 756/12.}

\maketitle

\begin{abstract}
We define Schwartz functions, tempered functions and tempered distributions on (possibly singular) real algebraic varieties. We prove that all classical properties of these spaces, defined previously on affine spaces and on Nash manifolds, also hold in the case of affine real algebraic varieties, and give partial results for the non-affine case.
\end{abstract}

\tableofcontents

\section{Introduction}

Schwartz functions are classically defined as smooth functions such that they, and all their (partial) derivatives, decay at infinity faster than the inverse of any polynomial. On $\bR$, for instance, a smooth function $f$ is called Schwartz if for any $n,k\in\mathbb{N}, x^n f^{(k)}$ is bounded (where $f^{(k)}$ is the $k^{th}$ derivative of $f$). This was formulated on $\bR^n$ by Laurent Schwartz, and later on Nash manifolds (smooth semi-algebraic varieties), see \cite{dC,AG}. As Schwartz functions are defined using algebraic notions, it is natural to define Schwartz spaces of real algebraic varieties -- this is the main goal of this paper.

We start by defining \emph{Schwartz functions} on a real algebraic set in $\bR^n$ as the quotient of the space of Schwartz functions on $\bR^n$ by the ideal of Schwartz functions which vanish identically on the set. As the space of Schwartz functions on $\bR^n$ is a Fr\'echet space, we consider the natural quotient topology. We show that this definition makes sense for an affine variety, i.e. it does not depend on the embedding. The independence of the embedding enables us to carry many properties of Schwartz functions shown in \cite{dC,AG} to the (singular) algebraic affine case.

A classical, very useful, characterization of Schwartz functions on a Nash open subset $U$ of a Nash manifold $X$ is $\mathcal{S}(U)\cong W_Z$. Here $\mathcal{S}(U)$ is the space of Schwartz functions on $U$, $Z=X\setminus U$ and $W_Z$ is the space of all Schwartz functions on $X$ which vanish identically with all their derivatives on $Z$. A smooth function that vanishes identically with all its derivatives at some point, is said to be \emph{flat} at this point. In this paper we prove a similar result, for $X$ being a real affine algebraic variety, and $U\subset X$ a Zariski open subset. In order to define $W_Z$ we have to make sense of the notion "derivative" at a singular point of a variety, or more specifically the notion of a "flat" function at such a point. We do this by the following (a-priori naive) definition: we say that a function $f$ on an algebraic set $X\subset\bR^n$ is flat at some $y\in X$, if it is the restriction of some $C^\infty(\bR^n)$ function, that is flat at $y$. This point of view suggests the characterization of Schwartz functions by {\bf local} means only, where the global conditions of "rapid decaying at infinity" are translated to local conditions of flatness at "all points added in infinity" added in some compactification process. We make this claim precise in Theorem \ref{char-schwartz-on-open-from-affine-complete} and Remark \ref{milman_comment}.

\emph{Tempered functions} on $\bR^n$ are classically defined as smooth functions
which are bounded by a polynomial, and all their (partial) derivatives are
bounded by polynomials as well. Similar definition works for Nash manifolds. As with Schwartz functions, we define tempered functions on an affine
algebraic variety as
restrictions of classical tempered functions on some embedding space. An
important tool we develop using this definition is the tempered partition
of unity.

We define \emph{tempered distributions} as the space of continuous linear functionals on the space of Schwartz functions (the dual space), and get an important consequence of the isomorphism $\mathcal{S}(U)\cong W_Z$ discussed above: the restriction of tempered distributions on an affine algebraic variety to an open subset is onto. This result is known for tempered distributions on Nash manifolds (see [AG]), but it does not hold for general distributions. This claim follows from the isomorphism discussed above and the Hahn-Banach Theorem. In fact, we show that in general tempered distributions form a flabby sheaf.

The main results appearing in this paper for {\bf affine varieties} are as follows:

\begin{enumerate}
\item Let $X\subset\mathbb{R}^n$ be an algebraic set, then $\mathcal{S}(X)$ is a Fr\'echet space (Lemma \ref{Schwartz-is-Frechet}).

\item Let $\varphi:X_1\to X_2$ be a biregular isomorphism between two algebraic sets $X_1\subset \bR^{n_1},X_2\subset \bR^{n_2}$. Then $\varphi^*|_{\mathcal{S}(X_2)}:\mathcal{S}(X_2)\to\mathcal{S}(X_1)$ is an isomorphism of Fr\'echet spaces. This implies that the definition of Schwartz functions on an affine algebraic variety does not depend on the embedding (Lemma \ref{prop-indep-in-embed}).

\item Tempered partition of unity: Let $X$ be an affine algebraic variety, and let $\{V_i\}_{i=1}^m$ be a Zariski open cover of $X$. Then, there exist tempered functions $\{\beta_i\}_{i=1}^m$ on $X$, such that $supp(\beta_i)\subset V_i$
and  $\sum\limits_{i=1}^m\beta_i=1$.
 Furthermore, for any $m$-tuple $(\beta_1,..,\beta_m)$ of tempered functions on $X$ satisfying these conditions, and for any $\varphi\in\mathcal{S}(X)$, $(\beta_i\cdot\varphi)|_{V_i}\in\mathcal{S}(V_i)$ (Proposition \ref{temp-part-unity-affine-case} and Corollary \ref{Cor-2-partition-uni-always-holds}).

\item For an affine algebraic variety $M$, and a Zariski closed subset of it $X\subset M$, the restriction $\phi\mapsto\phi|_X$  maps $\mathcal{S}(M)$ onto $\mathcal{S}(X)$ (Theorem \ref{rest.-from-closed}).

\item Let $X$ be an affine algebraic variety, and let $Z\subset X$ be some Zariski closed subset. Define $U:=X\setminus Z$ and $W_Z:=\{\phi\in\mathcal{S}(X)|\phi\text{ is flat on }Z\}$. Then $W_Z$ is a closed subspace of $\mathcal{S}(X)$ (and so it is a Fr\'echet space), and extension by zero $\mathcal{S}(U)\to W_Z$ is an isomorphism of Fr\'echet spaces, whose inverse is the restriction of functions (Theorem \ref{char-schwartz-on-open-from-affine-complete}). As a consequence, the restriction morphism of tempered distributions $\mathcal{S}^*(X)\to\mathcal{S}^*(U)$ is onto (Theorem \ref{theorem-chararterization-of-Schwartz-on-open-with-duals}).

\item Let $X$ be an affine algebraic variety. The assignment of the space of Schwartz functions (respectively tempered functions, tempered distributions) to any
open $U\subset X$, together with the extension by zero $Ext_U^V$ from $U$
to any other open $V\supset U$ (restriction of functions, restrictions of functionals from $\mathcal{S}^*(V)$ to $\mathcal{S}^*(U)$), form a flabby cosheaf (sheaf, flabby sheaf) on $X$ (Propositions \ref{Schwartz-is-a-cosheaf}, \ref{Tempered-is-a-sheaf}, and \ref{Temp-dist-is-a-sheaf}).
\end{enumerate}

For {\bf general (not necessarily affine) varieties}, we define Schwartz functions as sums of extensions by zero of Schwartz functions on affine subvarieties, and prove this definition is independent on the covering and that the space of Schwartz functions is a Fr\'echet space (Lemma \ref{schwartz-on-non-aff-is-well-def} and Definition \ref{def-schwartz-on-non-affine}). We show that restrictions of Schwartz functions on an algebraic variety to a closed subset are Schwartz function of this subset (Proposition \ref{rest.-from-closed-non-affine}). We also define the space of tempered functions (Lemma \ref{lemma-for-def-temp-on-non-affine} and Definition \ref{def-temp-on-non-affine}), and prove that the space of Schwartz
functions on an algebraic variety forms a module over the space of tempered
functions (Proposition \ref{temp.-times-Schwartz-is-Schwartz_non_affine}). We also prove that if $X$ is an algebraic variety, then the assignment of the space of tempered functions to any open $U\subset X$, together with the restriction of functions, form a sheaf on $X$ (Proposition \ref{Tempered-is-a-sheaf-for-non-affine}). We define the notion of a flat function at a point, and by doing so we are able to define $W_Z$ for non-affine varieties. Finally, we prove that if $X$ is an algebraic variety and $U\subset X$ is some Zariski open subset, then the extensions by zero of functions in $\mathcal{S}(U)$ lie in $W_Z$ (Proposition \ref{ext-by-zero-for-non-affine}).

{\bf Structure of this paper:} in {\bf Section \ref{prelim}} we present the preliminary definitions and results we will use in this paper, mainly from real algebraic geometry and Schwartz functions on Nash manifolds.

In {\bf Section \ref{affine-case-chapter}} we define the space of Schwartz functions on an affine algebraic variety, and study its properties. We start by showing that it is a Fr\'echet space and proving that a useful partition of unity exists. Afterwards
we define the notion of flat functions at a point
on an affine algebraic variety, and characterize the spaces of Schwartz functions
on Zariski open subsets of an affine algebraic variety. We also define tempered
distributions, and prove that the restriction morphism from the space of
tempered distributions on an affine algebraic variety to the space of tempered
distributions on an open subset of it is onto. The proofs of two key lemmas in Section \ref{affine-case-chapter}
require some tools from subanalytic geometry -- Appendix \ref{subanalytic-geometry-1}
is dedicated to presenting these tools and completing the two proofs.

As proving that tempered functions and tempered distributions form sheaves, and that Schwartz functions form a cosheaf is quite technical, a separated section is dedicated to these proofs -- this is {\bf Section \ref{section_sheaves}}.

In {\bf Section \ref{non-affine-case}} we define the spaces of Scwhartz functions and of tempered functions on an arbitrary (not necessarily affine) real algebraic variety, and repeat some of the results we proved in the affine case. We also shortly discuss the difficulty of generalizing the rest of these results to the non-affine case, and suggest an idea that might enable overcoming this difficulty.

Throughout this paper the base field is always $\bR$, and we always consider the Zariski topology, unless otherwise stated.

{\bf Acknowledgments.} First and foremost, we would like to thank our advisor, Prof. Dmitry Gourevitch, for presenting us the subject in question in this paper, and for suggesting us, with endless patience, which paths to take in order to answer it. We also thank Prof. William A. Casselman for suggesting the basic idea of defining Schwartz functions on algebraic sets as restrictions of Schwartz functions on the embedding space. Our gratitude also goes to Prof. Avraham Aizenbud for many useful ideas, and to Prof. Dmitry Novikov for listening to our ideas along the way and suggesting improvements. Finally, we would like to thank Prof. Pierre D. Milman for explaining how to apply his results in order to solve some of our problems, and for reading a preliminary version of this paper and suggesting the simplification of some definitions -- without his assistance,
this paper would not have been completed.

\section{Preliminaries}\label{prelim}

In this section we present the basic definitions and results we will use in this paper. These include some background in real algebraic geometry (\ref{prelim-alg-geom}), Schwartz functions and tempered functions on Nash manifolds (\ref{prelim-schwartz-on-nash}), and Fr\'echet spaces (\ref{prelim-frechet}).

\subsection{Real algebraic geometry}\label{prelim-alg-geom}

We start by recalling the basic definitions:

\subsubsection{Definitions (following \cite{BCR})}\label{regular-functions}
Let $X\subset\bR^n$ be an algebraic set (i.e. the zero locus of a family of polynomials in $\bR[x_1,..,x_n]$). Denote $I_{Alg}(X):=\{p\in\bR[x_1,...,x_n]:p|_X=0\}$. Define \emph{the coordinate ring of $X$} by $\bR[X]:=\bR[x_1,...,x_n]/I_{Alg}(X)$. Let $V$ be an open subset of $X$. A function $f:V\to \bR$ is called \emph{a regular function} if $f=\frac{g}{h}$, where $g,h\in\bR[X]$ and $h^{-1}(0)\cap V=\emptyset$. Note that the space of regular functions on $V$ forms a ring. Moreover, the assignment of such a ring to any open subset of $X$ defines a sheaf on $X$. We call this sheaf \emph{the sheaf of regular functions on $X$} and denote it by $\mathcal{R}_X$. A map $F:V\to \bR^m$ ($F(x)=(F_1(x),...,F_m(x))$) is called \emph{a regular map} if for any $1\leq i \leq m$: $F_i$ is a regular function. Let $Y\subset\bR^m$ be an algebraic set, and let $U$ be an open subset of $Y$. A map from $V$ to $U$ is called \emph{a biregular isomorphism} if it is a bijective regular map, whose inverse map is also regular. In that case we will say that $V$ is biregular isomorphic to $U$. An \emph{affine algebraic variety} is a topological space $X'$, equipped with a sheaf of real valued functions $\mathcal{R}_{X'}$, isomorphic (as a ringed space) to an algebraic set $X\subset \bR^n$ with its Zariski topology, equipped with its sheaf of regular functions $\mathcal{R}_X$. The sheaf $\mathcal{R}_{X'}$ is called the sheaf of regular functions on $X'$, and the topology of $X'$ is called the Zariski topology. An \emph{algebraic variety} is a topological space $X'$, equipped with a sheaf of real valued functions $\mathcal{R}_{X'}$, such that there exists a finite open cover $\{U_i\}_{i=1}^n$ of $X'$, with each $U_i$ equipped with the sheaf $\mathcal{R}_{X'}|_{U_i}$ being an affine algebraic variety. The sheaf $\mathcal{R}_{X'}$ is called the sheaf of regular functions on $X'$, and the topology of $X'$ is called the Zariski topology.

Remark. Note that unlike in the algebraically closed case, the ring of regular function on $\bR^n$ is not $\bR[x_1,...,x_n]$, e.g. $\frac{1}{x^2+1}$ is not a polynomial, but it is regular on $\bR$.

The following two propositions discuss the nature of algebraic subsets of $\bR^n$ and of their open subsets:

\subsubsection{Proposition (\cite[Proposition 2.1.3]{BCR})}\label{any-alg-set-is-a-zero-locus}
Let $X\subset\bR^n$ be an algebraic set. There exists $f\in\mathbb{R}[x_1,x_2,..,x_n]$ such that $X$ is the zero locus of $f$, i.e. $X=\{x\in\bR^n|f(x)=0\}$.

\subsubsection{Proposition (\cite[Proposition 3.2.10]{BCR})}\label{open-of-alg-is-aff}
Let $X\subset\bR^n$ be an algebraic set, and $U$ an open subset of $X$. Then $(U,\mathcal{R}_X|_U)$ is an affine algebraic variety (when we define for any open $U'\subset U\subset X$: $\mathcal{R}_X|_U(U'):=\mathcal{R}_X(U')$).

The following Proposition \ref{Zariski-is-Noetherian} is implicitly used in \cite{BCR} (see for instance \cite[Corollary 3.2.4]{BCR}). For the reader's convenience we give its detailed proof in Appendix \ref{proof-of-Notherianity}.

\subsubsection{Proposition (the Zariski topology is Noetherian)}\label{Zariski-is-Noetherian} Let $X$ be a real algebraic variety, $U\subset X$ a (Zariski) open subset, and $\{U_\alpha\}_{\alpha\in I}$ an open cover of $U$. Then there exists a finite subcover $\{U_{\alpha_i}\}_{i=1}^k$.

\subsubsection{Definition}\label{def-complete-affine-var} An affine algebraic variety is called \emph{complete} if any regular function on it is bounded.

\subsubsection{Remark}\label{remark-complete-from-bcr} Definition \ref{def-complete-affine-var} is a special case of \cite[Definition 3.4.10]{BCR}. Note that if $X$ is a complete affine algebraic variety then for any closed embedding $i:X\hookrightarrow \bR^n$, $i(X)$ is compact in the Euclidean topology on $\bR^n$.

\subsubsection{Proposition (Algebraic Alexandrov compactification -- \cite[Proposition 3.5.3]{BCR})}\label{alexandrov_compactification}
Let $X$ be an affine algebraic variety that is not complete, then there exists a pair $(\dot{X},i)$ such that:

\begin{enumerate}
 \item $\dot{X}$ is a complete affine algebraic variety.
 \item $i:X\to\dot{X}$ is a biregular isomorphism from $X$ onto $i(X)$.
 \item $\dot{X}\setminus i(X)$ consists of a single point.
\end{enumerate}

\subsection{Schwartz functions and tempered functions on Nash manifolds}\label{prelim-schwartz-on-nash}

The reader not interested in Nash manifolds may skip the following paragraph, and replace the terms "an affine Nash manifold" and "restricted topology" by "a smooth affine algebraic variety" and "Zariski topology" (respectively) everywhere.

A restricted topological space is a space satisfying all axioms of a topological space, where the demand of infinite unions of open subsets being open is replaced by the (weaker) demand of {\bf finite} unions of open subsets being open. Semi-algebraic subsets are subsets defined using finitely many polynomials equalities and inequalities. The collection of all Euclidean open semi-algebraic subsets of a semi-algebraic set forms a restricted topology. An affine Nash manifold is a ringed restricted topological space, isomorphic, as a ringed space, to a closed smooth semi-algebraic subset of $\bR^n$, with the sheaf of rings (in the restricted topology) of Nash functions (smooth semi-algebraic functions). A Nash manifold is a ringed restricted topological space that has a finite open cover, such that any component in this cover is an affine Nash manifold. The interested reader may find an extensive introduction is \cite{AG}. We will start with the basic definitions of Schwartz functions and tempered functions on affine Nash manifolds, following \cite{AG}:

\subsubsection{Definition (\cite[Definition 4.1.1]{AG})}\label{AG-Def-4.1.1}
Let $M$ be an affine Nash manifold. Define \emph{the space of Schwartz functions on $M$} by $$\mathcal{S}(M):=\{\phi\in C^\infty(M)|\text{for any Nash differential operator }D\text{ on }M,D\phi\text{ is bounded}\}.$$ Introduce a topology on this space by the following system of semi-norms: $$||\phi||_D:=\sup\limits_{x\in M}|D\phi(x)|.$$

\subsubsection{Definition (\cite[Definition 4.2.1]{AG})}\label{AG-Def-4.2.1}
Let $M$ be an affine Nash manifold. A function $\alpha\in C^\infty(M)$ is called \emph{tempered} if for any Nash differential operator $D$ on $M$, there exists a Nash function $f$ on $M$, such that $|D\alpha|\leq f$. The space of tempered functions on $M$ is denoted by $\mathcal{T}(M)$.

The following results will be of special importance for us:

\subsubsection{Proposition (\cite[Proposition 4.2.1]{AG})}\label{AG-Prop-4.2.1}
Let $M$ be an affine Nash manifold and $\alpha$ a tempered function on $M$. Then $\alpha\mathcal{S}(M)\subset\mathcal{S}(M)$.

\subsubsection{Theorem (Partition of unity -- \cite[Theorem 4.4.1]{AG})}\label{AG-Thm-4.4.1}
Let $M$ be an affine Nash manifold, and let $\{U_i\}_{i=1}^n$ be a finite open semi-algebraic cover of $M$. Then:

\begin{enumerate}
 \item there exist tempered functions $\alpha_1,\alpha_2,..,\alpha_n$ on $M$ such that $supp(\alpha_i)\subset U_i$ and $\sum\limits_{i=1}^n \alpha_i=1$.
 \item Moreover, $\alpha_i$ can be chosen in such a way that for any $\phi\in\mathcal{S}(M)$, $\alpha_i\cdot\phi\in\mathcal{S}(U_i)$.
\end{enumerate}

\subsubsection{Proposition (\cite[Proposition 4.5.3]{AG})}\label{AG-Prop-4.5.3}
Let $M$ be an affine Nash manifold. The assignment of the space of tempered functions on $U$, to any open semi-algebraic $U\subset M$, together with the usual restriction maps, defines a sheaf of algebras on $M$ (in the restricted topology).

\subsubsection{Theorem (\cite[Theorem 4.6.1]{AG})}\label{AG-Thm-4.6.1}
Let $M$ be an affine Nash manifold, and $Z\hookrightarrow M$ be a closed Nash submanifold. The restriction $\mathcal{S}(M)\to\mathcal{S}(Z)$ is defined, continuous and onto. Moreover, it has a section $s:\mathcal{S}(Z)\to\mathcal{S}(M)$ such that if $\phi\in \mathcal{S}(Z)$ is zero at some point $p$ with all its derivatives, then $s(\phi)$ is also zero at $p$ with all its derivatives.

\subsubsection{Theorem (\cite[Theorem 4.6.2]{AG})}\label{AG-Thm-4.6.2}
Let $M$ be an affine Nash manifold, and $Z\hookrightarrow M$ be a closed Nash submanifold. The restriction $\mathcal{T}(M)\to\mathcal{T}(Z)$ is defined, continuous and onto. Moreover, it has a section $s:\mathcal{T}(Z)\to\mathcal{T}(M)$ such that if $\alpha\in \mathcal{T}(Z)$ is zero at some point $p$ with all its derivatives, then $s(\alpha)$ is also zero at $p$ with all its derivatives.

\subsubsection{Theorem (Characterization of Schwartz functions on open subsets -- \cite[Theorem 5.4.1]{AG})}\label{AG-Thm-5.4.1}
Let $M$ be an affine Nash manifold, $Z\hookrightarrow M$ be a closed semi-algebraic subset, and $U=M\setminus Z$. Let $W_Z$ be the closed subspace of $\mathcal{S}(M)$ defined by $$W_Z:=\{\phi\in\mathcal{S}(M)|\phi\text{ vanishes with all its derivatives on }Z\}.$$ Then restriction and extension by zero give an isomorphism $\mathcal{S}(U)\cong W_Z$.

The last property we will strongly used is the following:

\subsubsection{Proposition (\cite[Corollary 4.1.2]{AG})}\label{AG-Cor-4.1.2}
Let $M$ be an affine Nash manifold. Then $\mathcal{S}(M)$ is a Fr\'echet space.

We will use this fact later on, but first give here a brief reminder on these spaces.

\subsection{Fr\'echet spaces}\label{prelim-frechet} A Fr\'echet space is a metrizable, complete locally convex topological vector space. It can be shown that the topology of a Fr\'echet space can always be defined by a countable family of semi-norms. We will use the following results:

\subsubsection{Proposition (cf. \cite[Chapter 10]{T})}\label{closed-sub-of-F-is-F}
A closed subspace of a Fr\'echet space is a Fr\'echet space (for the induced topology).

\subsubsection{Proposition (cf. \cite[Proposition 7.9 and Chapter 10]{T})}\label{quotient-of-F-is-F}
A quotient of a Fr\'echet space by a closed subspace is a Fr\'echet space (for the quotient topology). Moreover, let $F$ be a Fr\'echet space whose topology is defined by a basis of continuous semi-norms $\mathcal{P}$, $K\subset F$ a closed subspace, and $\phi:F\to F/K$ the canonical mapping of $F$ onto $F/K$. Then the topology on $F/K$ is defined by the basis of continuous semi-norms $\dot{p}(\dot x)=\inf\limits_{\phi(x)=\dot{x}}p(x)$, where $p\in\mathcal{P}$.

\subsubsection{Theorem (Banach open mapping -- \cite[Chapter 17, Corollary 1]{T})}\label{Banach-open-mapping}
A bijective continuous linear map from a Fr\'echet space to another Fr\'echet space is an isomorphism, i.e. its inverse is also continuous.

\subsubsection{Theorem (Hahn-Banach -- cf. \cite[Chapter 18]{T})}\label{Hahn-Banach}
Let $F$ be a Fr\'echet space, and $K\subset F$ a closed subspace. By Proposition \ref{closed-sub-of-F-is-F} $K$ is a Fr\'echet space (with the induced topology). Define $F^*$ (respectively $K^*$) to be the space of continuous linear functionals on $F$ (on $K$). Then the restriction map $F^*\to K^*$ is onto.

\section{The affine case}\label{affine-case-chapter}

\subsection{Definition}\label{definition-minus-one}
Let $X\subset\bR^n$ be an algebraic subset. Let $\mathcal{S}(\bR^n)$ be the space of classical real valued Schwartz functions on $\bR^n$, and $I_{Sch}(X)\subset\mathcal{S}(\bR^n)$ be the ideal of all Schwartz functions that vanish identically on $X$. Define \emph{the space of Schwartz functions on $X$} by $\mathcal{S}(X):=\mathcal{S}(\bR^n)/I_{Sch}(X)$, equipped with the quotient topology.

\subsubsection{Remark} Equivalently we may define $$\mathcal{S}(X):=\{f:X\to\bR:\exists \tilde f\in\mathcal{S}(\bR^n)\text{ such that }\tilde{f}|_{X}=f\},$$ but then the definition of the topology will be a bit more complicated: recall that the topology of $\mathcal{S}(\bR^n)$ is given by a system of semi-norms $|f|_D:=\sup\limits_{x\in\bR^n}|Df(x)|$, where $D$ is an algebraic differential operator on $\bR^n$. This enables us to construct a system of semi-norms on $\mathcal{S}(X)$ by $|f|_D:=\inf\{|\tilde f|_D\:|\:\tilde{f}\in\mathcal{S}(\bR^n)\:,\:\tilde{f}|_X=f\}$, where $D$ is an algebraic differential operator on $\bR^n$. By Proposition \ref{quotient-of-F-is-F}, the topology defined by this system of semi-norms is the same as the quotient topology given in Definition \ref{definition-minus-one}.

\subsection{Lemma}\label{Schwartz-is-Frechet} $\mathcal{S}(X)$ is a Fr\'echet space.

Proof: $\mathcal{S}(\bR^n)$ is a Fr\'echet space (see Proposition \ref{AG-Cor-4.1.2}), $I_{Sch}(X)$ is a closed subspace (as $I_{Sch}(X)=\bigcap\limits_{x\in X}\{f\in\mathcal{S}(\bR^n)|f(x)=0\}$, is an intersection of closed sets) and so the quotient is a Fr\'echet space (see Proposition \ref{quotient-of-F-is-F}). \qed

\subsection{Lemma}\label{rest.-from-open}
Let $X\subset\bR^n$ be an algebraic set, and $U\subset\bR^n$ be some open set containing $X$. Consider $\mathcal{S}(\bR^n)$ and $\mathcal{S}(U)$ as defined in \ref{AG-Def-4.1.1} ($U$ is an open semi-algebraic set, hence it may be considered as an affine Nash manifold). Let $I_{Sch}^U(X)\subset\mathcal{S}(U)$ be the ideal of all Schwartz functions on $U$ that vanish identically on $X$. Then $\mathcal{S}(X)\cong \mathcal{S}(U)/I_{Sch}^U(X)$ (isomorphism of Fr\'echet spaces).

Proof: By the same reasoning as in Lemma \ref{Schwartz-is-Frechet}, $\mathcal{S}(U)/I_{Sch}^U(X)$ is a Fr\'echet space. By Theorem \ref{AG-Thm-5.4.1}, $\mathcal{S}(U)$ is isomorphic to a closed subspace of $\mathcal{S}(\bR^n)$, and so by Proposition \ref{quotient-of-F-is-F} it is enough to check that  $\mathcal{S}(X):=\mathcal{S}(\bR^n)/I_{Sch}(X)$ and $\mathcal{S}(U)/I_{Sch}^U(X)$ are equal as sets, i.e. that a function on $X$ is a restriction of a Schwartz function on $\bR^n$ if and only if it is a restriction of a Schwartz function on $U$.

Let $f\in\mathcal{S}(U)|_X$. There exists $F\in\mathcal{S}(U)$ such that $F|_X=f$. By Theorem \ref{AG-Thm-5.4.1}, extending $F$ by zero to a function on $\bR^n$ (denote it by $\tilde{F}$) is a function in $\mathcal{S}(\bR^n)$. Then $f=\tilde{F}|_X$ and so $f\in\mathcal{S}(\bR^n)|_X$.

Let $f\in\mathcal{S}(\bR^n)|_X$. There exists $F\in\mathcal{S}(\bR^n)$ such that $F|_X=f$. Denote $U':=\bR^n\setminus X$. $\{U,U'\}$ form an open cover of $\bR^n$ and so, by Theorem \ref{AG-Thm-4.4.1}, there exist tempered functions $\alpha_1,\alpha_2$ such that $supp(\alpha_1)\subset U$, $supp(\alpha_2)\subset U'$ and $\alpha_1+\alpha_2=1$ as a real valued function on $\bR^n$.  Moreover, $\alpha_1$ and $\alpha_2$ can be chosen such that $(\alpha_1\cdot F)|_U\in \mathcal{S}(U)$. As $\alpha_1|_X=1$, it follows that $((\alpha_1\cdot F)|_U)|_X=(\alpha_1\cdot F)|_X=F|_X=f$, and so $f\in\mathcal{S}(U)|_X$. \qed

In a similar way we prove a version of Lemma \ref{rest.-from-open} for tempered functions:

\subsection{Lemma}\label{rest.-from-open-temp.}
Let $X\subset\bR^n$ be an algebraic set, and $U\subset\bR^n$ be some open set containing $X$. Consider $\mathcal{T}(\bR^n)$ and $\mathcal{T}(U)$ (the spaces of tempered functions on $\bR^n$ and on $U$, respectively) as defined in \ref{AG-Def-4.2.1} ($U$ is an open semi-algebraic set, hence it may be considered as an affine Nash manifold). Then a function $f:X\to\bR$ is a restriction of a function $F\in\mathcal{T}(\bR^n)$ if and only if it is a restriction of a function $\tilde F\in\mathcal{T}(U)$.

Proof: let $f:X\to\bR$ be a restriction of some function $F\in\mathcal{T}(\bR^n)$. By Proposition \ref{AG-Prop-4.5.3}, $F|_U\in\mathcal{T}(U)$, and clearly $f=(F|_U)|_X$, i.e. $f$ is a restriction of a tempered function on $U$.

Let $f:X\to\bR$ be a restriction of some function $F\in\mathcal{T}(U)$. Denote $U':=\bR^n\setminus X$. $\{U,U'\}$ form an open cover of $\bR^n$ and so, by Theorem \ref{AG-Thm-4.4.1}, there exist tempered functions $\alpha_1,\alpha_2\in\mathcal{T}(\bR^n)$ such that $supp(\alpha_1)\subset U$, $supp(\alpha_2)\subset U'$ and $\alpha_1+\alpha_2=1$ as a real valued function on $\bR^n$. As tempered functions on affine Nash manifolds form a sheaf (see Proposition \ref{AG-Prop-4.5.3}), $\alpha_1|_U\in\mathcal{T}(U)$ and as $\mathcal{T}(U)$ is an algebra, $\alpha_1|_U\cdot F\in\mathcal{T}(U)$. Moreover, defining $F':\bR^n\to\bR$ by $F'|_U:=F$ and $F'|_{\bR^n\setminus U}:=0$, then as $supp(\alpha_1)\subset U$, we get that $\alpha_1\cdot F'\in\mathcal{T}(\bR^n)$. Since $\alpha_1|_X=1$ we have $(\alpha_1\cdot F')|_X=f$, i.e. $f$ is a restriction of a tempered function on $\bR^n$. \qed

\subsection{Lemma} Let $\varphi:X_1\to X_2$ be a biregular isomorphism between two algebraic sets $X_1\subset\bR^{n_1}$ and $X_2\subset\bR^{n_2}$.

\subsubsection{}\label{prop-indep-in-embed}
$\varphi^*|_{\mathcal{S}(X_2)}:\mathcal{S}(X_2)\to\mathcal{S}(X_1)$ is an isomorphism of Fr\'echet spaces.

\subsubsection{}\label{temp-is-invariant} If $f:X_2\to \bR$ is a restriction of a tempered function on $\bR^{n_2}$ (see Definition \ref{AG-Def-4.2.1}) then $\varphi^*f:=f\circ \varphi$ is a restriction of a tempered function on $\bR^{n_1}$.

Proof: By definition for any $x\in X_1$, we have $\varphi(x)=(\frac{f_1(x)}{g_1(x)},..,\frac{f_{n_2}(x)}{g_{n_2}(x)})$, where $f_1,..,f_{n_2},g_1,..,g_{n_2}\in\bR[X_1]$ and for any $1\leq i\leq {n_2}$, $g_i^{-1}(0)\cap X_1=\emptyset$. By abuse of notation we choose some representatives in $\bR[x_1,..,x_{n_1}]$ and consider $f_1,..,f_{n_2},g_1,..,g_{n_2}$ as functions in $\bR[x_1,..,x_{n_1}]$. Define $U:=\{x\in\bR^{n_1}|\prod\limits_{i=1}^{n_2} g_i(x)\neq 0\}$. $U$ is open in $\bR^{n_1}$ (also in the Euclidean topology), $X_1$ is a closed subset of $U$, and $\varphi$ can be naturally extended to a regular map $\tilde{\varphi}:U\to\bR^{n_2}$ (by the same formula of $\varphi$). Note that $U$ is an affine Nash manifold.

Similarly to the construction of $U$ and $\tilde\varphi$ above, we may construct an open $V\subset \bR^{n_2}$ and a function $\phi:V\to\bR^{n_1}$ such that $\phi|_{X_2}=\varphi^{-1}$. Note that $\phi\neq \tilde\varphi^{-1}$: in general $\tilde\varphi$ is not a bijection and $U\not\cong V$.

Consider the following diagram, where $\alpha$ is defined by $\alpha(x,y):=(x,y+\tilde{\varphi}(x))$.

\begin{center}
\leavevmode
\xymatrix{
X_1 \ar@{^{(}->}[r] & U \ar@/^1pc/[rr]^{Id\times 0} \ar@/_1pc/[rr]_{Id\times \tilde{\varphi}} && U\times\bR^{n_2} \ar@(ul,ur)^{\alpha} }
\end{center}

Clearly $U\times \{0\}$ is an affine Nash manifold isomorphic to $U$. Denote $\tilde U:=\alpha(U\times\{0\})$, then $\alpha$ restricted to $U\times\{0\}$ is an isomorphism of the affine Nash manifolds $U\times\{0\}$ and $\tilde U$ --  the inverse map is given by $\alpha^{-1}(x,y):=(x,y-\tilde{\varphi}(x))$ . Thus we have: $$\mathcal{S}(X_1)\cong \mathcal{S}(U)/I_{Sch}^U(X_1)\cong\mathcal{S}(\tilde U)/I_{Sch}^{\tilde U}(\alpha(X_1\times\{0\}))=\mathcal{S}(\tilde U)/I_{Sch}^{\tilde U}((Id\times\varphi)(X_1)),$$ where the first equivalence is by Lemma \ref{rest.-from-open}, the second is due the fact that $U\cong U\times\{0\}\cong\tilde U$ and $\mathcal{S}(U)\cong\mathcal{S}(U\times\{0\})\cong\mathcal{S}(\tilde U)$, and the third follows from the fact that $\tilde \varphi|_{X_1}=\varphi$. As always $I_{Sch}^U(X)$ is the ideal in $\mathcal{S}(U)$ of Schwartz functions identically vanishing on $X$.

As $\tilde U$ is closed in $U\times \bR^{n_2}$ (as it is defined by polynomial equalities on $U\times \bR^{n_2}$), then by Theorem \ref{AG-Thm-4.6.1} and Proposition \ref{quotient-of-F-is-F} we get that $$\mathcal{S}(\tilde U)/I_{Sch}^{\tilde U}((Id\times\varphi)(X_1))\cong \mathcal{S}(U\times\bR^{n_2})/I_{Sch}^{U\times\bR^{n_2}}((Id\times\varphi)(X_1)).$$

Applying Lemma \ref{rest.-from-open} again for the open subset $U\times V \subset U\times\bR^{n_2}$ we get that $$\mathcal{S}(U\times\bR^{n_2})/I_{Sch}^{U\times\bR^{n_2}}((Id\times\varphi)(X_1))\cong \mathcal{S}(U\times V)/I_{Sch}^{U\times V}((Id\times\varphi)(X_1)),$$

and thus we obtain $$\mathcal{S}(X_1)\cong\mathcal{S}(U\times V)/I_{Sch}^{U\times V}((Id\times\varphi)(X_1)).$$

Repeating the above construction using the following diagram:

\begin{center}
\leavevmode
\xymatrix{
X_2 \ar@{^{(}->}[r] & V \ar@/^1pc/[rr]^{0\times Id} \ar@/_1pc/[rr]_{\phi\times Id} && \bR^{n_2}\times V }
\end{center}

yields:

$$\mathcal{S}(X_2)\cong\mathcal{S}(U\times V)/I_{Sch}^{U\times V}((\varphi^{-1}\times Id)(X_2)).$$

Clearly $(Id\times\varphi)(X_1)=(\varphi^{-1}\times Id)(X_2)$, and so $\mathcal{S}(X_1)\cong \mathcal{S}(X_2)$. Note that the isomorphism constructed is in fact the pull back by $\varphi$ from $\mathcal{S}(X_2)$ onto $\mathcal{S}(X_1)$. This proves \ref{prop-indep-in-embed}.

The proof of \ref{temp-is-invariant} is the same as the proof of \ref{prop-indep-in-embed}, where one should consider tempered functions instead of Schwartz functions, and use Lemma \ref{rest.-from-open-temp.} and Theorem \ref{AG-Thm-4.6.2} instead of Lemma \ref{rest.-from-open} and Theorem \ref{AG-Thm-4.6.1}. \qed

\subsection{Definition}\label{definition-zero}
Let $X$ be a real affine algebraic variety, and let $i:X\hookrightarrow\bR^n$ be a closed embedding. A function $f:X\to\bR$ is called \emph{a Schwartz function on $X$} if $i_*f:=f\circ i^{-1}\in\mathcal{S}(i(X))$. Denote the space of all Schwartz functions on $X$ by $\mathcal{S}(X)$, and define a topology on $\mathcal{S}(X)$ by declaring a subset $U\subset\mathcal{S}(X)$ to be open if $i_*(U)\subset\mathcal{S}(i(X))$ is an open subset. By Lemma \ref{prop-indep-in-embed}, $\mathcal{S}(X)$ is well defined (independent of the embedding chosen).
\subsubsection{Remarks}
\begin{enumerate}
 \item If $X\cong\bR^m$ then Definition \ref{definition-zero} coincides with the classical one.
 \item If $X$ is smooth then Definition \ref{definition-zero} coincides with Definition \ref{AG-Def-4.1.1} of Schwartz functions on Nash manifolds (by substituting $M=\bR^n$ and $Z=i(X)$ in Theorem \ref{AG-Thm-4.6.1}).
\end{enumerate}

\subsection{Theorem}\label{rest.-from-closed}
Let $M$ be an affine algebraic variety, and let $X\subset M$ be a closed subset. Then the restriction from $M$ to $X$ defines an isomorphism $\mathcal{S}(X)\cong\mathcal{S}(M)/I_{Sch}^M(X)$ (with the quotient topology), where $I_{Sch}^M(X)$ is the ideal in $\mathcal{S}(M)$ of functions identically vanishing on $X$.

Proof: Take some closed embedding $M\hookrightarrow \bR^n$, then $X\hookrightarrow M\hookrightarrow \bR^n$ are closed embeddings. Then $\mathcal{S}(M)/I_{Sch}^M(X)=(\mathcal{S}(\bR^n)/I_{Sch}(M))/I_{Sch}^M(X)\cong\mathcal{S}(\bR^n)/I_{Sch}(X)=\mathcal{S}(X)$. \qed

\subsubsection{Remark} In particular for any $\phi\in\mathcal{S}(M)$, one has that $\phi|_X \in\mathcal{S}(X)$, and this restriction map $\mathcal{S}(M)\to\mathcal{S}(X)$ is onto.

\subsection{Definition}\label{temp-func-def} Let $X$ be an affine algebraic variety. A function $f:X\to\bR$ is called \emph{a tempered function on $X$}, if there exists a closed embedding $i:X\hookrightarrow \bR^n$ such that $i_* f:=f\circ i^{-1}$ is a restriction of a tempered function on $\bR^n$ to $i(X)$. By Lemma \ref{temp-is-invariant} in that case this property holds for any closed embedding. The set of all tempered functions form a unitary algebra, which we denote by $\mathcal{T}(X)$.

\subsection{Proposition}\label{temp.-times-Schwartz-is-Schwartz} Let $X$ be an affine algebraic variety, $t\in\mathcal{T}(X)$ and $s\in\mathcal{S}(X)$. Then $t\cdot s\in\mathcal{S}(X)$.

Proof: Consider some closed embedding $i:X\hookrightarrow\bR^n$ and identify $i(X)$ with $X$ (by definitions of tempered and Schwartz functions the choice of the embedding does not matter). There exist $T\in\mathcal{T}(\bR^n)$ and $S\in\mathcal{S}(\bR^n)$ such that $t=T|_X$ and $s=S|_X$. By Proposition \ref{AG-Prop-4.2.1}, $T\cdot S\in\mathcal{S}(\bR^n)$, and so $(T\cdot S)|_X=t\cdot s\in\mathcal{S}(X)$. \qed

\subsection{Corollary from Proposition \ref{open-of-alg-is-aff}}\label{cor.-defined-on-open}
Let $X$ be an affine algebraic variety, and $U$ an open subset of $X$. By Proposition \ref{open-of-alg-is-aff}, $(U,\mathcal{R}_X|_U)$ is an affine algebraic variety, and we may define $\mathcal{S}(U)$.

\subsection{Proposition (tempered partition of unity)}\label{temp-part-unity-affine-case}
Let $X$ be an affine algebraic variety, and let $\{V_i\}_{i=1}^{m}$ be a finite open cover of $X$. Then:
\begin{enumerate}
 \item There exist tempered functions $\{\beta_i\}_{i=1}^{m}$ on $X$, such that $supp(\beta_i)\subset V_i$ \footnote{A-priori the notion of support of a function on an affine algebraic variety $X$ is not well defined -- we consider the support defined by {\bf some Euclidean topology} on $X$, where we choose some closed embedding of $X$ in an affine space. It can be easily shown the this support is {\bf independent} of the embedding chosen, and so this notion is well defined. For more details see Remark \ref{remark-on-supp}.} and $\sum\limits_{i=1}\limits^{m}\beta_i=1$.
 \item We can choose $\{\beta_i\}_{i=1}^{m}$ in such a way that for any $\varphi\in\mathcal{S}(X)$, ($\beta_i\varphi)|_{V_i}\in\mathcal{S}(V_i)$.
\end{enumerate}

Proof: Consider some closed embedding $X\hookrightarrow \bR^n$. For any $1\leq i\leq m$ let $U_i\subset\bR^n$ be some open subset such that $V_i=X\cap U_i$. Define $U_{m+1}:=\bR^n\setminus X$ we get that $\{U_i\}_{i=1}^{m+1}$ is an open cover of $\bR^n$ (by Nash submanifolds). By Theorem \ref{AG-Thm-4.4.1}, there exist $\{\alpha_i\}_{i=1}^{m+1}$, tempered functions on $\bR^n$, such that $supp(\alpha_i)\subset U_i$, $\sum\limits_{i=1}\limits^{m+1}\alpha_i=1$, and $\{\alpha_i\}_{i=1}^{m+1}$ can be chosen in such a way that for any $\psi\in\mathcal{S}(\bR^n)$, ($\alpha_i\psi)|_{U_i}\in\mathcal{S}(U_i)$. For $1\leq i \leq m+1$ define $\beta_i:=\alpha_i|_{X}$. Clearly for $1\leq i \leq m$, $supp(\beta_i)\subset V_i$. Since $\alpha_{m+1}|_X=0$ then $\beta_{m+1}|_X=0$ and so $\sum\limits_{i=1}\limits^{m}\beta_i=1$. By definition $\{\beta_i\}_{i=1}^{m}$ are tempered functions on $X$. This proves (1).

Now consider $\varphi\in\mathcal{S}(X)$. By definition there exists $\tilde{\varphi}\in\mathcal{S}(\bR^n)$ such that $\varphi=\tilde{\varphi}|_X$ and for $1\leq i\leq m$, $(\alpha_i \tilde{\varphi})|_{U_i}\in\mathcal{S}(U_i)$. By Theorem \ref{rest.-from-closed}, as $V_i$ is closed in $U_i$, we get that $((\alpha_i \tilde{\varphi})|_{U_i})|_{V_i}\in\mathcal{S}(V_i)$. But $((\alpha_i \tilde{\varphi})|_{U_i})|_{V_i}=(\beta_i \varphi)|_{V_i}$, and so (2) is proven. \qed

\subsection{Definition}\label{flat_new_def} Let $X\subset\bR^n$ be an algebraic set and $y\in X$ be some point. A function $f:X\to \bR$ is  \emph{flat at $y$} if there exists $F\in C^\infty(\bR^n)$ with $f=F|_X$, such that the Taylor series of $F$ at $y$ is identically zero. If $f$ is flat at $y$ for any $y\in Z$ (where $Z\subset X$ is some subset),
we will say that $f$ is \emph{flat at $Z$}.

\subsubsection{An important remark} As the Taylor series is only dependent on the Euclidean local behaviour
of functions, one may substitute $\bR^n$ above by any Euclidean open subset of $\bR^n$
containing $X$. This will be done in \ref{appendix-res-from-open-nbrhd},
and will be used when it will be more convenient.

\subsubsection{Warning} $f$ is flat at $Z$ means that $f$ is flat at $y$ for any $y\in Z$. It does not mean, a-priori, that there exists $F\in C^\infty(\bR^n)$ with $f=F|_X$,
such that all the Taylor series of $F$, at any point $y\in
Z$ are identically
zero. Lemma \ref{Milman_main_lemma} addresses this matter.

The proofs of the following Lemmas \ref{Milman_main_lemma} and \ref{X-flat-is-invariant} are given in Appendix \ref{subanalytic-geometry-2},
as they require some tools from subanalytic geometry.

\subsection{Lemma}\label{Milman_main_lemma}Let $X$ be a compact (in the
Euclidean topology) algebraic set in $\bR^n$, and let $Z\subset X$ be some
(Zariski) closed subset. Define $U:=X\setminus Z$, $$W_Z:=\{\phi:X\to\bR|\exists\tilde\phi\in
C^\infty(\bR^n)\text{ such that }\tilde\phi|_{X}=\phi\text{ and }\phi\text{
is flat on }Z\}$$ and $$(W^{\bR^n}_Z)^{comp}:=\{\phi\in C^\infty(\bR^n)|\phi\text{
is a compactly supported and is flat on }Z\}.$$ Then, for any $f\in W_Z$, there
exists $\tilde{f}\in (W_{Z}^{\bR^n})^{comp}$ such that $\tilde{f}|_{X}=f$.

\subsection{Lemma}\label{X-flat-is-invariant} Let $\varphi:X_1\to X_2$ be a biregular isomorphism between two algebraic sets $X_1\subset\bR^{n_1}$ and $X_2\subset\bR^{n_2}$. If $f:X_2\to \bR$ is flat at some $p\in X_2$ then $\varphi^*f:=f\circ \varphi$ is flat at $\varphi^{-1}(p)$.

\subsection{Definition} Let $X$ be an affine algebraic variety, and let $f:X\to\bR$ be some function. We say that \emph{$f$ is flat at $p\in X$} if there exists a closed embedding $i:X\hookrightarrow \bR^n$ such that $i_* f:=f\circ i^{-1}:i(X)\to\bR$ is flat at $i(p)$. By Lemma \ref{X-flat-is-invariant} in that case this property holds for any closed embedding.

\subsection{Proposition (extension by zero)}\label{ext-by-zero}
Let $X$ be an affine algebraic variety, and $U$ an open subset of $X$. By Corollary \ref{cor.-defined-on-open}, $\mathcal{S}(U)$ is defined. Then the extension by zero to $X$ of a Schwartz function on $U$ is a Schwartz function on $X$, which is flat at $X\setminus U$.

Proof: $X$ is affine, thus we may choose some closed embedding $X\hookrightarrow \bR^n$, and so we may think of $X$ as an algebraic set. According to Proposition \ref{any-alg-set-is-a-zero-locus}, there exists $F\in \bR[x_1,...,x_n]$ such that $X$ is the zero locus of $F$ (denote $X=zeros(F)$). $U\subset X$ is Zariski open in $X$, i.e. $Z:=X\setminus U$ is Zariski closed in $X$, thus $Z$ is Zariski closed in $\bR^n$. As before, there exists $G\in \bR[x_1,...,x_n]$ such that $Z=zeros(G)$. Define $V:=\bR^n\setminus Z$. Note that $U=X\setminus Z$ is a closed subset of $V$, as $X\setminus Z=X\cap V$. As $V$ is open in $\bR^n$, by Proposition \ref{open-of-alg-is-aff}, $V$ is an affine variety. Consider some closed embedding $V\hookrightarrow \bR^m$. By Theorem \ref{rest.-from-closed} $\mathcal{S}(U)\cong \mathcal{S}(V)|_U$. Let $h\in \mathcal{S}(U)$, then there exists $\bar{h}\in \mathcal{S}(V)$ such that $h=\bar{h}|_U$. As $V$ is an open subset of the Nash manifold $\bR^n$, by Theorem \ref{AG-Thm-5.4.1}, the extension of $\bar{h}$ by zero to $\bR^n$ (denote it by $\hat{h}$) is a Schwartz function on $\bR^n$ which is flat on $Z$. Finally, defining $\tilde{h}:=\hat{h}|_X$, we get that $\tilde{h}\in\mathcal{S}(X)$ (by Definition \ref{definition-zero}), $\tilde{h}|_U=h$ (by definition), and $\tilde{h}$ is flat on $X\setminus U$ (as $\hat{h}$ is an extension of $\tilde{h}$ to $\bR^n$ which is flat at $X\setminus U$). \qed

\subsection{Lemma}\label{lemma-shcwartz-is-like-smooth} Let $X$ be a compact (in the Euclidean topology) algebraic set in $\bR^n$, then $$\mathcal{S}(X)=\{f:X\to\bR:\exists \tilde f\in C^\infty(\bR^n)\text{ such that }\tilde{f}|_{X}=f\}.$$

Proof: The inclusion $\subset$ is trivial, as $\mathcal{S}(\bR^n)\subset C^\infty(\bR^n)$. For the inclusion $\supset$ take some $g\in\{f:X\to\bR:\exists \tilde f\in C^\infty(\bR^n)\text{ such that }\tilde{f}|_{X}=f\}$. Denote by $\tilde g$ some $C^\infty(\bR^n)$ function satisfying $\tilde g|_X=g$. Let $\rho\in C^\infty(\bR^n)$ be a compactly supported (in the Euclidean topology) function such that $\rho|_X=1$ (it is standard to show such $\rho$ exists by convolving the characteristic function of some bounded open subset containing $X$ with some appropriate approximation of unity). Then $\rho\cdot \tilde g$ is a smooth compactly supported function on $\bR^n$, hence $\rho\cdot \tilde g\in\mathcal{S}(\bR^n)$. Moreover $(\rho\cdot\tilde g)|_X=\tilde g|_X=g$, and so $g\in\mathcal{S}(X)$. \qed

\subsection{Lemma}\label{Milman_main_lemma_tag} Let $X$ be a compact (in the Euclidean topology) algebraic set in $\bR^n$, and let $Z\subset X$ be some (Zariski) closed subset. Define $U:=X\setminus Z$, $W_Z:=\{\phi\in\mathcal{S}(X)|\phi\text{
is flat on }Z\}$ and $W^{\bR^n}_Z:=\{\phi\in\mathcal{S}(\bR^n)|\phi\text{
is flat on }Z\}$. Then for any $f\in W_Z$, there
exists $\tilde{f}\in W_{Z}^{\bR^n}$ such that $\tilde{f}|_{X}=f$.

Proof: This is immediate from Lemma \ref{Milman_main_lemma} and Lemma \ref{lemma-shcwartz-is-like-smooth}. \qed

\subsection{Proposition}\label{rest.-from-W_z}
Let $X$ be an affine algebraic variety, and let $Z\subset X$ be some closed subset. Define $U:=X\setminus Z$ and $W_Z:=\{\phi\in\mathcal{S}(X)|\phi\text{ is flat on }Z\}$. Then restriction from $X$ to $U$ of a function in $W_Z$ is a Schwartz function on $U$, i.e. $Res_{X}^{U}(W_Z)\subset\mathcal{S}(U)$.

Proof: We first prove the case where $X$ is complete, and then deduce the non-complete case from the complete case:

Consider $X$ as an algebraic subset in $\bR^n$.

Case 1 -- $X$ is complete. By Remark \ref{remark-complete-from-bcr}, $X$ is compact in the Euclidean topology in $\bR^n$. Define $U^{\bR^n}:=\bR^n\setminus Z$ and $W^{\bR^n}_Z:=\{\phi\in\mathcal{S}(\bR^n)|\phi\text{ is flat on }Z\}$. As $Z$ is closed in $X$, which is closed in $\bR^n$, $Z$ is also closed in $\bR^n$, and so $U^{\bR^n}$ is open in $\bR^n$. As $U=U^{\bR^n}\cap X$, we get that $U$ is closed in $U^{\bR^n}$. We will show that indeed $Res_{X}^{U}(W_Z)\subset\mathcal{S}(U)$ by showing the existence of the following 3 maps:

\begin{center}
\leavevmode
\xymatrix{
& W^{\bR^n}_Z \ar[dr]^{Res_{\bR^n}^{U^{\bR^n}}}_{(2)} \ar@{>>}[dl]^{(1)}_{Res_{\bR^n}^{X}} \\
W_Z \ar@{.>}[dr]_{Res_{X}^{U}}  & & \mathcal{S}(U^{\bR^n}) \ar[ld]^{Res_{U^{\bR^n}}^{U}}_{(3)} \\
& \mathcal{S}(U)}
\end{center}

\

Clearly a restriction of a function in $W_{Z}^{\bR^n}$ to $X$ lies in $W_Z$, i.e. $Res_{\bR^n}^{X}$ is well defined. By Lemma \ref{Milman_main_lemma_tag}, $Res_{\bR^n}^{X}$ is onto. Let $g\in W_{Z}^{\bR^n}$. Then, by Theorem \ref{AG-Thm-5.4.1}, $g|_{U^{\bR^n}}\in\mathcal{S}(U^{\bR^n})$, i.e. map (2) is well defined. Let $h\in\mathcal{S}(U^{\bR^n})$. Then, by Theorem \ref{rest.-from-closed}, $h|_U\in\mathcal{S}(U)$, i.e. map (3) is well defined. Thus Proposition \ref{rest.-from-W_z} holds if $X$ is complete.

Case 2 -- $X$ is non-complete. Consider a one point compactification, i.e. a pair $(\dot{X},i)$ as in Proposition \ref{alexandrov_compactification},
and take some $f\in W_Z\subset \mathcal{S}(X)$. As $i:X\to i(X)$ is a biregular
isomorphism, $i_*f:=f\circ i^{-1}\in\mathcal{S}(i(X))$. As $i(X)$ is open
in $\dot X$, by Proposition \ref{ext-by-zero}
there exists $\dot{f}\in\mathcal{S}(\dot{X})$ such that $i_*f=\dot{f}|_{i(X)}$
($\dot{f}$ is the extension by zero to $\dot X$ of $i_*f$).

Denote $p:=\dot{X}\setminus i(X)$ (by Proposition \ref{alexandrov_compactification},
$p$ is a single point, hence it is closed in $\dot X$). We claim that $i(Z)\cup
\{p\}$ is closed in $\dot{X}$. Indeed, $\dot X\setminus (i(Z)\cup \{p\})=i(X)\setminus
i(Z)$ is open in $i(X)$ (as $Z$ is closed in $X$ and $i$ is a biregular isomorphism
of $X$ and $i(X)$), and $i(X)$ is open in $\dot X$. Now define $U':=\dot
X \setminus (i(Z)\cup \{p\})$, is open in $\dot X$. By case 1,
$Res_{\dot{X}}^{U'}(\dot{f})\in \mathcal{S}(U')$. Observe that $i^{-1}|_{U'}$
is a biregular isomorphism of $U'$ and $U$, and so $(i^{-1}|_{U'})_*Res_{\dot{X}}^{U'}(\dot{f})\in
\mathcal{S}(U)$. But $(i^{-1}|_{U'})_*Res_{\dot{X}}^{U'}(\dot{f})=(i^{-1}|_{U'})_*((i_*f)|_{U'})=f|_U$,
thus $f|_U\in\mathcal{S}(U)$. \qed

\subsection{Theorem}\label{char-schwartz-on-open-from-affine-complete}
Let $X$ be an affine algebraic variety, and let $Z\subset X$ be some closed subset. Define $U:=X\setminus Z$ and $W_Z:=\{\phi\in\mathcal{S}(X)|\phi\text{ is flat on }Z\}$. Then $W_Z$ is a closed subspace of $\mathcal{S}(X)$ (and so it is a Fr\'echet space), and extension by zero $\mathcal{S}(U)\to W_Z$ (denote $Ext_{U}^{X}:\mathcal{S}(U)\to W_Z$) is an isomorphism of Fr\'echet spaces, whose inverse is the restriction of functions (denote $Res_{X}^{U}:W_Z\to\mathcal{S}(U)$).

Proof: $W_Z$ is clearly a subspace of $\mathcal{S}(X)$. As $W_z=\bigcap\limits_{z\in Z}\{\phi\in\mathcal{S}(X)|\phi\text{ is flat on }z\}$ is an intersection of closed sets, it is closed. Thus $W_Z$ is a Fr\'echet space (see Proposition \ref{closed-sub-of-F-is-F}).

By Proposition \ref{ext-by-zero} the extension of a function in $\mathcal{S}(U)$ by zero to $X$ is a function in $\mathcal{S}(X)$ that is $X$-flat at $Z$, i.e. $Ext_U^X(\mathcal{S}(U))\subset W_Z$. Furthermore, we claim that $Ext_U^X$ is continuous.

Indeed, take some closed embedding $X\hookrightarrow \bR^n$ and consider $X$ as an algebraic set. Define $W=\bR^n\setminus Z$. $W$ is an affine Nash manifold containing $U$ (which is an affine algebraic variety), and $U=W\cap X$, i.e. $U$ is closed in $W$. Now take some closed embedding $W\hookrightarrow \bR^N$, then by Theorem \ref{rest.-from-closed} and Proposition \ref{quotient-of-F-is-F}, $\mathcal{S}(U)\cong \mathcal{S}(W)/I_{Sch}^{W}(U)$. As $W$ is open in $\bR^n$, by Theorem \ref{AG-Thm-5.4.1}, $Ext_{W}^{\bR^n}$ is a closed embedding $\mathcal{S}(W)\hookrightarrow \mathcal{S}(\bR^n)$ (and in particular it is a continuous map). Then we may write $I_{Sch}^{W}(U)=I_{Sch}^{\bR^n}(X)\cap \mathcal{S}(W)$. In particular $I_{Sch}^{W}(U)$ is closed in $\mathcal{S}(\bR^n)$. Finally, as the embedding $\mathcal{S}(W)\hookrightarrow \mathcal{S}(\bR^n)$ is continuous, then the map $$\mathcal{S}(U)=\mathcal{S}(W)/(I_{Sch}^{\bR^n}(X)\cap \mathcal{S}(W))\to \mathcal{S}(\bR^n)/I_{Sch}^{\bR^n}(X)=\mathcal{S}(X)$$ is continuous as well, i.e. $Ext_U^X$ is continuous.

By Proposition \ref{rest.-from-W_z} the restriction of a function in $W_Z$ to $U$ is a Schwartz function on $U$, i.e. $Res_{X}^{U}(W_Z)\subset\mathcal{S}(U)$.

By definition, $Res_X^U\circ Ext_U^X:\mathcal{S}(U)\to\mathcal{S}(U)$ is the identity operator on $\mathcal{S}(U)$ and $Ext_U^X\circ Res_X^U:W_Z\to W_Z$ is the identity operator on $W_Z$. Thus, $Ext_U^X$ is a continuous (linear) bijection. Then, by Banach open mapping Theorem (\ref{Banach-open-mapping}), $Ext_U^X$ is an isomorphism of Fr\'echet spaces. \qed

\subsection{Corollary}\label{cor-rami-def-is-equiv} Let $X$ be an affine algebraic variety. A Schwartz function $f\in\mathcal{S}(X)$ is flat at $p\in X$ if and only if $f|_{X\setminus \{p\}}\in \mathcal{S}({X\setminus \{p\}})$.

Proof: Apply Theorem \ref{char-schwartz-on-open-from-affine-complete} to $Z=\{p\}$. \qed

\subsubsection{Remark} By the same argument for an arbitrary function $f\in C^\infty(X)$ (i.e. a function that is a restriction of a smooth function from an open neighborhood of some closed embedding of $X$) and any $p\in X$, the following conditions are equivalent:

\begin{enumerate}
 \item $f$ is flat at $p$.
 \item There exists a smooth compactly supported function $\rho$ on some affine space which $X$ is closely embedded in, such that $\rho$ is identically 1 on some open neighborhood of $p$ and $(f\cdot\rho)|_{X\setminus \{p\}}\in \mathcal{S}({X\setminus \{p\}})$.
 \item For any smooth compactly supported function $\rho$ on any affine space which $X$ is closely embedded in, such that $\rho$ is identically 1 on some open neighborhood of $p$, one has $(f\cdot\rho)|_{X\setminus \{p\}}\in \mathcal{S}({X\setminus \{p\}})$.
\end{enumerate}

Theorem \ref{char-schwartz-on-open-from-affine-complete} also implies that part (2) of Proposition \ref{temp-part-unity-affine-case} holds for any partition of unity:

\subsection{Corollary}\label{Cor-2-partition-uni-always-holds} Let $X$ be an affine algebraic variety, and let $V\subset X$ be some open subset of $X$. Then, for any $\beta\in\mathcal{T}(X)$ such that $supp(\beta)\subset V$, and for any $\varphi\in\mathcal{S}(X)$, one has $(\beta\cdot\varphi)|_{V}\in\mathcal{S}(V)$.

Proof: By Proposition \ref{temp.-times-Schwartz-is-Schwartz} $\beta\cdot\varphi\in\mathcal{S}(X)$. Consider $X$ as an algebraic subset of some $\bR^n$ (i.e. choose some closed embedding). By definition $supp(\beta)$ is a closed subset of $\bR^n$ in the Euclidean topology. There exists some (Zariski) open $\tilde V\subset \bR^n$ such that $V=\tilde V\cap X$. As $supp(\beta)\subset \tilde V$, which is also an open subset of $\bR^n$ in the Euclidean topology, the function $\beta\cdot\varphi$ is flat on $X\setminus V$. Thus, by Theorem \ref{char-schwartz-on-open-from-affine-complete}, $(\beta\cdot\varphi)|_{V}\in\mathcal{S}(V)$. \qed

\subsection{Remark}\label{milman_comment} Theorem \ref{char-schwartz-on-open-from-affine-complete} suggests the following point of view on Schwartz functions: given an affine algebraic variety $X$ we take some affine compactification of it, i.e. we consider $X$ as an open subset of some complete affine variety $Y$ \footnote{We used above Alexandrov (one point) compactification, but this is not necessary.}. Then, a Schwartz function on $X$ is just a smooth function on $Y$ (in the sense it is the restriction to $Y$ of a smooth function on the ambient space of $Y$), that is flat on $Y\setminus X$. This point of view is convenient, as it involves only local properties -- the condition of "rapidly decaying at infinity" is translated to the condition of flatness at "all points added in infinity" in the compactification process. This is also true in the Nash category (by Theorem \ref{AG-Thm-5.4.1}), and the easiest example is the case where $X=\bR$, where one can identify $\bR$ with the unit circle without a point.

\subsection{Definition}\label{def-temp-dist}
Let $X$ be an affine algebraic variety. Define the \emph{space of tempered distributions on $X$} as the space of continuous linear functionals on $\mathcal{S}(X)$. Denote this space by $\mathcal{S}^*(X)$.

\subsection{Theorem}\label{theorem-chararterization-of-Schwartz-on-open-with-duals}
Let $X$ be an affine algebraic variety, and let $U\subset X$ be some (Zariski) open subset. Then extending any function in $\mathcal{S}(U)$ by zero on $X\setminus U$ is a closed embedding $\mathcal{S}(U)\hookrightarrow \mathcal{S}(X)$, and the restriction morphism $\mathcal{S}^*(X)\to\mathcal{S}^*(U)$ is onto.

Proof: The first part of the Theorem is just a restatement of Theorem \ref{char-schwartz-on-open-from-affine-complete} (substituting $Z=X\setminus U$). The second part follows from the fact that $\mathcal{S}(X)$ is a Fr\'echet space (\ref{Schwartz-is-Frechet}) and from the Hahn-Banach Theorem (\ref{Hahn-Banach}). \qed

\section{Sheaf and cosheaf properties}\label{section_sheaves}

This section is devoted to proving that tempered functions and tempered distributions form sheaves (Propositions \ref{Tempered-is-a-sheaf} and \ref{Temp-dist-is-a-sheaf}) and that Schwartz functions form a cosheaf (Proposition \ref{Schwartz-is-a-cosheaf}). The precise definition of cosheaves is given in Appendix \ref{sheaves-and-cosheaves}.

\subsection{Lemma (restrictions of tempered functions to closed and to open subsets)}\label{temp.-rest} Let $X$ be an affine algebraic variety, and let $U\subset X$ be some open subset. Then $Res_X^{X\setminus U}(\mathcal{T}(X))=\mathcal{T}(X\setminus U)$ and $Res_X^U(\mathcal{T}(X))\subset\mathcal{T}(U)$.

Proof: Consider some closed embedding of $X$ in some affine space, i.e. consider $X$ as an algebraic subset of $\bR^n$. Then $\mathcal{T}(X):=Res_{\bR^n}^X(\mathcal{T}(\bR^n))$. In these settings $X\setminus U$ is also an algebraic subset of $\bR^n$, and so $\mathcal{T}(X\setminus U):=Res_{\bR^n}^{X\setminus U}(\mathcal{T}(\bR^n))=Res^{X\setminus U}_{X}(Res_{\bR^n}^X(\mathcal{T}(\bR^n)))=Res_X^{X\setminus U}(\mathcal{T}(X))$. This proves the first part of Lemma \ref{temp.-rest}.

As $U\subset X$ is open, there exists an open subset $\tilde U\subset \bR^n$ such that $U=\tilde U\cap X$, and $U$ is closed in $\tilde U$. In particular $\tilde U$ is a Nash submanifold of $\bR^n$. Let $t\in\mathcal{T}(X)$. By definition there exists $T\in\mathcal{T}(\bR^n)$ such that $T|_X=t$. By Proposition \ref{AG-Prop-4.5.3} $T|_{\tilde U}\in \mathcal{T}(\tilde U)$. As $\tilde U$ is open in $\bR^n$ it is also an affine algebraic variety (by Proposition \ref{open-of-alg-is-aff}). So $\tilde U$ is an affine algebraic variety and $U$ is a closed subset of it, thus by the first part of the lemma $t|_U=(T|_{\tilde U})|_U\in\mathcal{T}(U)$. \qed

\subsection{Corollary}\label{can-restrict-tempered} Let $X$ be an affine algebraic variety, and let $U,V\subset X$ be two open subsets of $X$ such that $U\subset V$. Then $Res_V^U(\mathcal{T}(V))\subset \mathcal{T}(U)$.

Proof: By Proposition \ref{open-of-alg-is-aff} $V$ is an affine algebraic variety, so we may apply Lemma \ref{temp.-rest} to the affine algebraic variety $V$ and its open subset $U$. \qed

\subsection{Proposition}\label{Tempered-is-a-sheaf} Let $X$ be an affine algebraic variety. The assignment of the space of tempered functions to any open $U\subset X$, together with the restriction of functions, form a sheaf on $X$.

Proof: By Corollary \ref{can-restrict-tempered} the above is a pre-sheaf.

Let $\{U_i\}_{i\in I}$ be some open cover of $U$, and $t,s\in\mathcal{T}(U)$ be such that for any $i\in I$, $t|_{U_i}=s|_{U_i}$. As $t$ and $s$ are simply real valued functions on $U$, clearly $s=t$, and so the axiom of uniqueness holds.

Now let $t_i\in\mathcal{T}(U_i)$ be such that for any $i,j\in I$, $t_i|_{U_i\cap U_j}=t_j|_{U_i\cap U_j}$. Clearly there exists a (unique) function $t:U\to\bR$ such that for any $i\in I$, $t|_{U_i}=t_i$. In order to prove that the existence axiom holds, it is thus left to show that $t\in\mathcal{T}(U)$. As always we may consider $X$ as an algebraic subset of $\bR^n$. By Proposition \ref{Zariski-is-Noetherian} we may assume $|I|<\infty$ by choosing some subcover and showing $t|_{U_i}=t_i$ only for indices $i$ in this subcover (as the functions we begin with agree on the intersections, this will automatically hold for all the other indices we omitted). By standard induction on the number of indices (i.e. the number of sets in the chosen finite subcover), it is enough to show that the following holds:

Let $X\subset \bR^n$ be an algebraic subset and let $U_1,U_2\subset X$ be two open subsets. Assume that for any $i\in\{1,2\}$ we are given $t_i\in \mathcal{T}(U_i)$ such that $t_1|_{U_1\cap U_2}=t_2|_{U_1\cap U_2}$. Then, there exists a function $t\in\mathcal{T}(U_1\cup U_2)$ such that $t|_{U_1}=t_1,t|_{U_2}=t_2$.

Clearly there exists a (unique) function $t:U_1\cup U_2\to\bR$ such that $t|_{U_i}=t_i$. It is left to show that $t\in\mathcal{T}(U_1\cup U_2)$. Indeed, there exist open sets $\tilde U_i\subset \bR^n$ such that $U_i=\tilde U_i\cap X$, and $U_i$ is closed in $\tilde U_i$. Then, by Lemma \ref{rest.-from-open-temp.}, there exist $T_i\in\mathcal{T}(\tilde U_i)$ such that $t_i=T_i|_{U_i}$. Define $U=U_1\cup U_2$ and $\tilde U=\tilde U_1\cup \tilde U_2$. As $\tilde U$ is an affine Nash manifold, and $\{\tilde U_1, \tilde U_2\}$ is an open cover of $\tilde U$, by Theorem \ref{AG-Thm-4.4.1} there exist $\alpha_1,\alpha_2\in\mathcal{T}(\tilde U)$ such that $supp(\alpha_i)\subset \tilde U_i$ and $\alpha_1+\alpha_2=1$. Define $T'_1(x):=T_1(x)$ if $x\in\tilde U_1$ and $T'_1(x):=0$ if $x\in \tilde U\setminus\tilde U_1$. Similarly define $T'_2$ on $\tilde U$ by extending $T_2$ by zero. Define a new function on $\tilde U$ by $T:=\alpha_1\cdot T'_1+\alpha_2\cdot T'_2$.

By Proposition \ref{AG-Prop-4.5.3}, in order to show that $T\in\mathcal{T}(\tilde U)$, it is enough to show that for $i\in\{1,2\}$, $T|_{\tilde U_i}\in\mathcal{T}(\tilde U_i)$. Let us show this for $i=1$ (symmetrical arguments work for $i=2$): $T|_{\tilde U_1}=\alpha_1|_{\tilde U_1}\cdot T'_1|_{\tilde U_1}+(\alpha_2\cdot T'_2)|_{\tilde U_1}=\alpha_1|_{\tilde U_1}\cdot T_1+(\alpha_2\cdot T'_2)|_{\tilde U_1}$. As the space of tempered functions is an algebra, it is enough to show each of these three functions belongs to $\mathcal{T}(\tilde U_1)$. By Proposition \ref{AG-Prop-4.5.3} $\alpha_1|_{\tilde U_1}\in\mathcal{T}(\tilde U_1)$. By construction $T_1\in\mathcal{T}(\tilde U_1)$. In order to show that $(\alpha_2\cdot T'_2)|_{\tilde U_1}\in\mathcal{T}(\tilde U_1)$ we use Proposition \ref{AG-Prop-4.5.3} again: as $\{\tilde U_1\cap \tilde U_2,\tilde U_1\setminus(supp(\alpha_2)\cap \tilde U_1)\}$ is an open cover of $\tilde U_1$, it is enough to show that $((\alpha_2\cdot T'_2)|_{\tilde U_1})|_{\tilde U_1\cap \tilde U_2}\in\mathcal{T}(\tilde U_1\cap \tilde U_2)$ and $((\alpha_2\cdot T'_2)|_{\tilde U_1})|_{\tilde U_1\setminus(supp(\alpha_2)\cap \tilde U_1)}\in\mathcal{T}(\tilde U_1\setminus(supp(\alpha_2)\cap \tilde U_1))$. The later is obvious as $((\alpha_2\cdot T'_2)|_{\tilde U_1})|_{\tilde U_1\setminus(supp(\alpha_2)\cap \tilde U_1)}=0$, and the first also holds as both $\alpha_2|_{\tilde U_1\cap\tilde U_2}\in \mathcal{T}(\tilde U_1\cap\tilde U_2)$ and $T_2|_{\tilde U_1\cap\tilde U_2}\in \mathcal{T}(\tilde U_1\cap\tilde U_2)$ as $\tilde U_1\cap\tilde U_2$ is open in both $\tilde U$ and in $\tilde U_2$ (and again by Proposition \ref{AG-Prop-4.5.3}). Finally, as $T\in\mathcal{T}(\tilde U)$ and $U_1\cup U_2=U\subset \tilde U$ is a closed subset, then by Lemma \ref{temp.-rest}, $t=T|_{U_1\cup U_2}\in\mathcal{T}(U_1\cup U_2)$. \qed

\subsection{Proposition}\label{Temp-dist-is-a-sheaf} Let $X$ be an affine algebraic variety. The assignment of the space of tempered distributions to any open $U\subset X$, together with restrictions of functionals from $\mathcal{S}^*(U)$ to $\mathcal{S}^*(V)$, for any other open $V\subset U$, form a sheaf on $X$.

Proof: Any open $U\subset X$ is an affine algebraic variety (by Proposition \ref{open-of-alg-is-aff}), and any open $V\subset X$ contained in $U$ is open in $U$. Thus, by Proposition \ref{ext-by-zero} the above is a pre-sheaf.

Let $\{U_i\}_{i\in I}$ be some open cover of $U$. By Proposition \ref{Zariski-is-Noetherian} there exists a finite open cover $\{U_i\}_{i=1}^k$. Note that $U$ is an affine algebraic variety (by Proposition \ref{open-of-alg-is-aff}). Then, by Proposition \ref{temp-part-unity-affine-case}, for any $1\leq i\leq k$ there exists $\beta_i\in\mathcal{T}(U)$ such that $supp(\beta_i)\subset U_i$, $\sum\limits_{i=1}^k \beta_i=1$ and for any $s\in\mathcal{S}(U)$, $(\beta_i\cdot s)|_{U_i}\in\mathcal{S}(U_i)$.

Now let $\xi,\zeta\in\mathcal{S}^*(U)$ be such that for any $i\in I$, $\xi|_{\mathcal{S}(U_i)}=\zeta|_{\mathcal{S}(U_i)}$. In particular, as $U_i$ is open in $U$ (and as $\{1,2,..,k\}\subset I$), then for any $1\leq i\leq k$, $\xi|_{\mathcal{S}(U_i)}=\zeta|_{\mathcal{S}(U_i)}$. Let $s\in \mathcal{S}(U)$. Note that $s=\sum\limits_{i=1}^{k} (\beta_i\cdot s)$, so we may calculate: $$\xi(s)-\zeta(s)=\xi(\sum\limits_{i=1}^{k} (\beta_i\cdot s))-\zeta(\sum\limits_{i=1}^{k} (\beta_i\cdot s))=\sum\limits_{i=1}^k(\xi(\beta_i\cdot s)-\zeta(\beta_i\cdot s))=0,$$  where the second equality follows from linearity of $\xi$ and $\zeta$, and the third equality follows from the facts that $(\beta_i\cdot s)|_{U_i}\in\mathcal{S}(U_i)$ and $\xi|_{\mathcal{S}(U_i)}=\zeta|_{\mathcal{S}(U_i)}$. Thus we showed that the axiom of uniqueness holds.

Now let $\xi_i\in\mathcal{S}^*(U_i)$ be such that for any $i,j\in I$, $\xi_i|_{\mathcal{S}(U_i\cap U_j)}=\xi_j|_{\mathcal{S}(U_i\cap U_j)}$. In particular, as $U_i\cap U_j$ is open in $U$ (and as $\{1,2,..,k\}\subset I$), then for any $1\leq i\leq k$, $\xi_i|_{\mathcal{S}(U_i\cap U_j)}=\xi_j|_{\mathcal{S}(U_i\cap U_j)}$. We define a functional $\xi\in\mathcal{S}^*(U)$ by the following formula for any $s\in\mathcal{S}(U)$: $$\xi(s)=\xi(\sum\limits_{i=1}^{k} (\beta_i\cdot s)):=\sum\limits_{i=1}^k \xi_i(\beta_i\cdot s).$$

We claim that for any $\alpha\in I$, one has $\xi|_{\mathcal{S}(U_\alpha)}=\xi_\alpha$. Indeed, $\{U_\alpha\cap U_i\}_{i=1}^k$ is an open cover of the affine algebraic variety $U_\alpha$. Note that $\{\beta_i|_{U_\alpha}\}_{i=1}^k$ is "a partition of unity"  of $U_\alpha$ as defined in Proposition \ref{temp-part-unity-affine-case}, i.e. $\beta_i|_{U_\alpha}\in\mathcal{T}(U_\alpha)$ (this follows from Proposition \ref{Tempered-is-a-sheaf}), $\sum\limits_{i=1}^k \beta_i|_{U_\alpha}=1$, $supp(\beta_i|_{U_\alpha})\subset U_\alpha\cap U_i$ and (by Corollary \ref{Cor-2-partition-uni-always-holds}) for any $s\in\mathcal{S}(U_\alpha)$ one has $(\beta_i|_{U_\alpha}\cdot s)|_{U_\alpha\cap U_i}\in\mathcal{S}(U_\alpha\cap U_i)$. Also note that for any $1\leq i\leq k$, one has $\xi_s|_{U_\alpha \cap U_i}=\xi_i|_{U_\alpha \cap U_i}$. Finally, we are ready to calculate (for any $s\in\mathcal{S}(U_\alpha)$, where we also think of $s$ as a function in $\mathcal{S}(U)$, by the usual extension by zero): $$\xi_\alpha(s)=\xi_\alpha(\sum\limits_{i=1}^k \beta_i\cdot s)=\sum\limits_{i=1}^k \xi_\alpha(\beta_i\cdot s)=\sum\limits_{i=1}^k \xi_i(\beta_i\cdot s)=\xi(s),$$ i.e. the axiom of existence holds. \qed

\subsection{Proposition}\label{Schwartz-is-a-cosheaf} Let $X$ be an affine algebraic variety. The assignment of the space of Schwartz functions to any open $U\subset X$, together with the extension by zero $Ext_U^V$ from $U$ to any other open $V\supset U$, form a cosheaf on $X$.

Proof: Any open $V\subset X$ is an affine algebraic variety (by Proposition \ref{open-of-alg-is-aff}), and any open $U\subset X$ contained in $V$ is open in $V$. Thus, by Theorem \ref{theorem-chararterization-of-Schwartz-on-open-with-duals} the above is a pre-cosheaf.

Let $\{U_i\}_{i\in I}$ be some open cover of $U$, and let $s\in\mathcal{S}(U)$. By Proposition \ref{Zariski-is-Noetherian} there exists a finite subcover $\{U_i\}_{i=1}^k$. By Proposition \ref{open-of-alg-is-aff} and Corollary \ref{cor.-defined-on-open} we may apply Proposition \ref{temp-part-unity-affine-case} on $U$ and so there exist $\beta_1,..\beta_k\in\mathcal{T}(U)$ such that $\sum\limits_{i=1}^k\beta_i=1$, and for any $1\leq i \leq k$, $supp(\beta_i)\subset U_i$, and $(\beta_i\cdot s)|_{U_i}\in\mathcal{S}(U_i)$. Then we may write: $$s=\sum\limits_{i=1}^k \beta_i\cdot s=\sum\limits_{i=1}^k Ext_{U_i}^{U}((\beta_i\cdot s)|_{U_i}),$$ and so axiom (1) holds.

Now assume we are given $s_i\in\mathcal{S}(U_i)$ for any $i\in J\subset I$, where $1<|J|<\infty$, such that $\sum\limits_{i\in J} Ext_{U_i}^U (s_i)=0$. We want to prove that for any $i>j\in J$ there exists $s_{i,j}\in\mathcal{S}(U_i\cap U_j)$ such that for any $i\in J$, $s_i=\sum\limits_{i>j\in J}Ext_{U_i\cap U_j}^{U_i}(s_{i,j})-\sum\limits_{i<j\in J}Ext_{U_i\cap U_j}^{U_i}(s_{j,i})$, and so axiom (2) holds. We prove this claim by induction on $|J|$. For $|J|=2$ one has $s_1|_{U_1\cap U_2}=-s_2|_{U_1\cap U_2}$, so defining $s_{2,1}=s_2|_{U_1\cap U_2}$ the claim holds -- the only non trivial fact to verify is that $s_2|_{U_1\cap U_2}\in\mathcal{S}(U_1\cap U_2)$: indeed, by Proposition \ref{ext-by-zero}, $-Ext_{U_1}^{U_1\cup U_2}(s_1)$ is a Schwartz function on $U_1\cup U_2$ that is flat on $U_1\cup U_2\setminus U_1$, and $Ext_{U_2}^{U_1\cup U_2}(s_2)$ is a Schwartz function on $U_1\cup U_2$ that is flat on $U_1\cup U_2\setminus U_2$. But as $Ext_{U_2}^{U_1\cup U_2}(s_2)=-Ext_{U_1}^{U_1\cup U_2}(s_1)$, we have that $Ext_{U_2}^{U_1\cup U_2}(s_2)$ is flat on $(U_1\cup U_2)\setminus (U_1\cap U_2)$. Then, by Theorem \ref{char-schwartz-on-open-from-affine-complete}, $s_2|_{U_1\cap U_2}=(Ext_{U_2}^{U_1\cup U_2}(s_2))|_{U_1\cap U_2}$ is a Schwartz function on $U_1\cap U_2$.

Now assume the claim holds for any $J$ of cardinality up to $k$, and let $J=\{1,2,..,k,k+1\}$. Without loss of generality we may assume $U=\bigcup\limits_{i=1}^{k+1} U_i$, and so we have for any $1\leq i\leq k+1$, $s_i\in\mathcal{S}(U_i)$ such that $\sum\limits_{i=1}^{k+1} Ext_{U_i}^U (s_i)=0$. Define $\tilde U:=\bigcup\limits_{i=1}^k U_i$. Note that $s_{k+1}|_{U_{k+1}\setminus(U_{k+1}\cap\tilde U)}=0$. As $\{U_i\}_{i=1}^{k}$ is an open cover of the affine $\tilde U$, then by Proposition \ref{temp-part-unity-affine-case} there exist $\{\beta_i\}_{i=1}^{k}\subset \mathcal{T}(\tilde U)$ such that for any $1\leq i\leq k$, $supp(\beta_i)\subset U_i$ and $\sum\limits_{i=1}^k \beta_i=1$.

Let $x\in U_{k+1}\setminus (U_{k+1}\cap \tilde U)$. By Proposition \ref{ext-by-zero}, for any $1\leq i\leq k$, $Ext_{U_i}^{U}(s_i)$ is flat at $x$. Then $Ext_{U_{k+1}}^{U}(s_{k+1})=-\sum\limits_{i=1}^{k}Ext_{U_i}^{U}(s_i)$ is also flat at $x$. Applying Theorem \ref{char-schwartz-on-open-from-affine-complete} (note that $U\setminus(U_{k+1}\setminus(\tilde U\cap U_{k+1}))=\tilde U$) we get that $(Ext_{U_{k+1}}^{U}(s_{k+1}))|_{\tilde U}\in\mathcal{S}(\tilde U)$. Now by Corollary \ref{Cor-2-partition-uni-always-holds} we have: $$(\beta_i|_{U_i}\cdot Ext_{U_i\cap U_{k+1}}^{U_i}(s_{k+1}|_{U_i\cap U_{k+1}}))=\beta_i|_{ U_{i}}\cdot (Ext_{U_{k+1}}^{\tilde U}(s_{k+1}))|_{ U_{i}}\in\mathcal{S}(U_{i}).$$

Define for any $1\leq i\leq k$, $$\gamma_i:=s_i+(\beta_i|_{U_i}\cdot Ext_{U_i\cap U_{k+1}}^{U_i}(s_{k+1}|_{U_i\cap U_{k+1}})).$$ Note that $\gamma_i\in\mathcal{S}(U_i)$ and that $\sum\limits_{i=1}^kExt_{U_i}^{\tilde U}(\gamma_i)=0$. Thus, by induction hypothesis, for any $1\leq j<i\leq k$ there exist $s_{i,j}\in\mathcal{S}(U_i\cap U_j)$ such that for any $1\leq i\leq k $, $$\gamma_i=\sum\limits_{i>j\geq 1}Ext_{U_i\cap U_j}^{U_i}(s_{i,j})-\sum\limits_{i<j\leq k}Ext_{U_i\cap U_j}^{U_i}(s_{j,i}).$$ For any $1\leq i\leq k$, define $s_{k+1,i}:=\beta_i|_{ U_{k+1}\cap U_i}\cdot s_{k+1}|_{ U_{k+1}\cap U_i}$. Then $\gamma_i=s_i+Ext_{U_{k+1}\cap U_i}^{U_i}(s_{{k+1},i})$, where $U_{k+1}\cap U_i$ is open in $U_i$. As both $\gamma_i$ and $s_i$ lie in $\mathcal{S}(U_i)$, so does $Ext_{U_{k+1}\cap U_i}^{U_i}(s_{{k+1},i})=\gamma_i-s_i$.

We claim that $s_{k+1,i}\in\mathcal{S}(U_{k+1}\cap U_i)$: denoting $f:=\beta_i|_{U_i}\cdot Ext_{U_i\cap U_{k+1}}^{U_i}(s_{k+1}|_{U_i\cap U_{k+1}})$, we saw above that $f\in\mathcal{S}(U_i)$. Thus, as $s_{k+1,i}=f|_{U_i\cap U_{k+1}}$, by Theorem \ref{char-schwartz-on-open-from-affine-complete} we need to show that $f$ is flat at $U_i\setminus (U_{k+1}\cap U_i)=(U_{k+1}\cup U_i)\setminus U_{k+1}$. Define $g:=Ext_{U_{k+1}}^{U_{k+1}\cup U_i}(s_{k+1})$. Then, by Theorem \ref{char-schwartz-on-open-from-affine-complete}, $g\in\mathcal{S}(U_i\cup U_{k+1})$, and $g$ is flat at $(U_i\cup U_{k+1})\setminus U_{k+1}$. In particular, $\tilde g:=g|_{U_i}$ is flat at $(U_i\cup U_{k+1})\setminus U_{k+1}$. Let $x\in (U_i\cup U_{k+1})\setminus U_{k+1}$, and let $\rho\in\mathcal{S}(U_i)$ be "a bump function around $x$", i.e. a restriction to $U_i$ of a smooth compactly supported function on some affine space which $U_i$ is closely embedded in, such that $\rho=1$ on some Euclidean open neighborhood of $x$. Then, by Corollary \ref{cor-rami-def-is-equiv}, $(\rho\cdot\tilde g)|_{U_i\setminus \{x\}}\in\mathcal{S}(U_i\setminus \{x\})$. By Proposition \ref{Tempered-is-a-sheaf}, $\beta_i|_{U_i\setminus \{x\}}\in\mathcal{T}(U_i\setminus \{x\})$. Thus, by Proposition \ref{temp.-times-Schwartz-is-Schwartz}, $(\beta_i|_{U_i}\cdot\rho\cdot\tilde g)|_{U_i\setminus \{x\}}\in\mathcal{S}(U_i\setminus \{x\})$. On the one hand, by Theorem \ref{char-schwartz-on-open-from-affine-complete}, $Ext_{U_i\setminus\{x\}}^{U_i}((\beta_i|_{U_i}\cdot\rho\cdot\tilde g)|_{U_i\setminus \{x\}})\in\mathcal{S}(U_i)$. On the other hand $(\beta_i|_{U_i}\cdot\rho\cdot\tilde g)|_{U_i}$ is a continuous function on $U_i$ that equals to $Ext_{U_i\setminus\{x\}}^{U_i}((\beta_i|_{U_i}\cdot\rho\cdot\tilde g)|_{U_i\setminus \{x\}})$ on $U_i\setminus \{x\}$. We deduce that $Ext_{U_i\setminus\{x\}}^{U_i}((\beta_i|_{U_i}\cdot\rho\cdot\tilde g)|_{U_i\setminus \{x\}})=(\beta_i|_{U_i}\cdot\rho\cdot\tilde g)|_{U_i}$, and so $(\beta_i|_{U_i}\cdot\rho\cdot\tilde g)|_{U_i}$ is flat at $x$. Finally, as flatness in a Euclidean local property, and as $\rho$ equals 1 on some Euclidean neighborhood of $x$, it follows that $\beta_i|_{U_i}\cdot\tilde g=f$ is flat at $x$.

Then it is easily seen that for any $1\leq i\leq k$ we have $$s_i=\sum\limits_{i>j\in \{1,2,..,k+1\}}Ext_{U_i\cap U_j}^{U_i}(s_{i,j})-\sum\limits_{i<j\in \{1,2,..,k+1\}}Ext_{U_i\cap U_j}^{U_i}(s_{j,i}).$$ It is left to check that $s_{k+1}=\sum\limits_{i=1}^{k}Ext_{U_{k+1}\cap U_i}^{U_{k+1}}(s_{k+1,i})$. Indeed, $$\sum\limits_{i=1}^{k}Ext_{U_{k+1}\cap U_i}^{U_{k+1}}(s_{k+1,i})=\sum\limits_{i=1}^{k}Ext_{U_{k+1}\cap U_i}^{U_{k+1}}(\beta_i|_{ U_{k+1}\cap U_i}\cdot s_{k+1}|_{ U_{k+1}\cap U_i})=$$ $$s_{k+1}\cdot\sum\limits_{i=1}^{k}Ext_{U_{k+1}\cap U_i}^{U_{k+1}}(\beta_i|_{ U_{k+1}\cap U_i})=s_{k+1}.$$ \qed

\section{The general case}\label{non-affine-case}

For any (non necessarily affine) algebraic variety $X$ denote the space of all real valued functions on $X$ by $Func(X,\bR)$.

\subsection{Lemma}\label{schwartz-on-non-aff-is-well-def}Let $X$ be an algebraic variety, and let $X=\bigcup\limits_{i=1}^kX_i=\bigcup\limits_{i=k+1}^lX_i$ be two open covers of $X$ by affine algebraic varieties. There is a natural map $\phi_1:\bigoplus\limits_{i=1}^kFunc(X_i,\bR)\to Func(X,\bR)$ and a natural map $\phi_2:\bigoplus\limits_{i=k+1}^lFunc(X_i,\bR)\to Func(X,\bR)$. Then $\phi_1(\bigoplus\limits_{i=1}^k\mathcal{S}(X_i))\cong\bigoplus\limits_{i=1}^k\mathcal{S}(X_i)/Ker(\phi_1|_{\bigoplus\limits_{i=1}^k\mathcal{S}(X_i)})$ has a natural structure of a Fr\'echet space, and there is an isomorphism of Fr\'echet spaces $\phi_1(\bigoplus\limits_{i=1}^k\mathcal{S}(X_i))\cong\phi_2(\bigoplus\limits_{i=k+1}^l\mathcal{S}(X_i))$.

Proof: It follows from Proposition \ref{quotient-of-F-is-F} that $\phi_1(\bigoplus\limits_{i=1}^k\mathcal{S}(X_i))\cong(\bigoplus\limits_{i=1}^k\mathcal{S}(X_i))/Ker(\phi_1|_{\bigoplus\limits_{i=1}^k\mathcal{S}(X_i)})$ is indeed a Fr\'echet space: a direct sum of Fr\'echet spaces is clearly a Fr\'echet space, and the kernel of $\phi_1|_{\bigoplus\limits_{i=1}^k\mathcal{S}(X_i)}$ is a closed subspace, as $\bigoplus\limits_{i=1}^k s_i\in Ker(\phi_1|_{\bigoplus\limits_{i=1}^k\mathcal{S}(X_i)})$ if and only if for any $x\in X$, $\sum\limits_{i\in J_x}s_i(x)=0$, where $J_x:=\{1\leq i\leq k|x\in X_i\}$, i.e. the kernel is given by infinitely many "closed conditions".

Note that $X=\bigcup\limits_{i=1}^k\bigcup\limits_{j=k+1}^lX_i\cap X_j$ is an open cover of $X$ by affine algebraic varieties. There is a natural map $\phi_3:\bigoplus\limits_{i=1}^k\bigoplus\limits_{j=k+1}^lFunc(X_i\cap X_j,\bR)\to Func(X,\bR)$. It is therefore enough to prove that $\phi_1(\bigoplus\limits_{i=1}^k\mathcal{S}(X_i))\cong\phi_3(\bigoplus\limits_{i=1}^k\bigoplus\limits_{j=k+1}^l\mathcal{S}(X_i\cap X_j))$. Note that $\{X_i\cap X_j\}_{j=k+1}^{l}$ is an open cover of $X_i$, and denote the natural map $\phi^i:\bigoplus\limits_{j=k+1}^{l}Func(X_i\cap X_j,\bR)\to Func(X_i,\bR)$. As $\phi_3=\phi_1\circ\bigoplus\limits_{i=1}^{k}\phi^i$, it is enough to prove that $\bigoplus\limits_{j=k+1}^{l}\mathcal{S}(X_i\cap X_j)/Ker(\phi^i|_{\bigoplus\limits_{j=k+1}^{l}\mathcal{S}(X_i\cap X_j)})\cong\mathcal{S}(X_i)$. Indeed, first we note that we have an equality of sets by Proposition \ref{Schwartz-is-a-cosheaf} and its proof: the  inclusion $\supset$ follows from axiom (1), and the inclusion $\subset$ follows from the fact that $Ext_{X_i\cap X_j}^{X_i}(\mathcal{S}(X_i\cap X_j))\subset\mathcal{S}(X_i)$. By Theorem \ref{char-schwartz-on-open-from-affine-complete} this extension is a closed embedding $\mathcal{S}(X_i\cap X_j)\hookrightarrow\mathcal{S}(X_i)$, and in particular it is continuous. Thus we have, by Theorem \ref{Banach-open-mapping}, an isomorphism of Fr\'echet spaces. \qed

\subsection{Definition}\label{def-schwartz-on-non-affine} Let $X$ be an algebraic variety, and let $X=\bigcup\limits_{i=1}^kX_i$ be some open cover of $X$ by affine algebraic varieties. There is a natural map $\phi:\bigoplus\limits_{i=1}^kFunc(X_i,\bR)\to Func(X,\bR)$. Define \emph{the space of Schwartz functions on $X$} by $\mathcal{S}(X):=(\bigoplus\limits_{i=1}^k\mathcal{S}(X_i))/Ker(\phi|_{\bigoplus\limits_{i=1}^k\mathcal{S}(X_i)})$, with the natural quotient topology. By Lemma \ref{schwartz-on-non-aff-is-well-def} this definition is independent of the cover chosen, and $\mathcal{S}(X)$ is a Fr\'echet space.

\subsection{Proposition}\label{rest.-from-closed-non-affine}Let $X$ be an algebraic variety, and $Z\subset X$ be some Zariski closed subset. Then $Res_{X}^{Z}(\mathcal{S}(X))\subset\mathcal{S}(Z)$.

Proof: Let $s\in\mathcal{S}(X)$, and let $X=\bigcup\limits_{i=1}^{n}X_i$ be some affine cover, such that
$s=\sum\limits_{i=1}^{n}Ext_{X_i}^{X}(s_i)$ for some $s_i\in\mathcal{S}(X_i)$. $Z\cap X_i$ is open in $Z$ and closed in $X_i$ (thus affine), and so $Z=\bigcup\limits_{i=1}^{n}Z\cap X_i$ is an affine open cover. By Theorem \ref{rest.-from-closed} for any $1\leq i\leq n$,  $s_i|_{Z\cap X_i}\in\mathcal{S}(Z\cap X_i)$, and thus  $s|_Z=\sum_{i=1}^{n}Ext_{X_i\cap Z}^{Z}(s|_{X_i\cap Z})\in\mathcal{S}(Z)$. \qed

\subsection{Lemma}\label{lemma-for-def-temp-on-non-affine} Let $X$ be an algebraic variety, and let $t:X\to\bR$ be some function. Then the following conditions are equivalent:

\begin{enumerate}
\item There exists an open affine cover $X=\bigcup\limits_{i=1}^{k}X_i$ such that for any $1\leq i\leq k$, $t|_{X_i}\in\mathcal{T}(X_i)$.
\item For any open affine cover $X=\bigcup\limits_{i=1}^{k}X_i$ and any $1\leq i\leq k$, $t|_{X_i}\in\mathcal{T}(X_i)$.
\end{enumerate}

Proof: Clearly (2) implies (1). For the other side assume there exist two open affine covers $X=\bigcup\limits_{i=1}^{k}X_i=\bigcup\limits_{j=k+1}^{l}X_i$ such that for any $k+1\leq j\leq l$, $t|_{X_j}\in\mathcal{T}(X_j)$. Fix some $1\leq i\leq k$. Note that $\{X_i\cap X_j\}_{j=k+1}^{l}$ is an open cover of $X_i$. By Proposition \ref{Tempered-is-a-sheaf}, as for any $k+1\leq j\leq l$, $t|_{X_j}\in\mathcal{T}(X_j)$, we have $t|_{X_j\cap X_i}\in\mathcal{T}(X_j\cap X_i)$. Applying Proposition \ref{Tempered-is-a-sheaf} once again, we get that $t|_{X_i}\in\mathcal{T}(X_i)$. \qed

\subsection{Definition}\label{def-temp-on-non-affine} Let $X$ be an algebraic variety. A real valued function $t:X\to\bR$ is called \emph{a tempered function on $X$} if it satisfies the equivalent conditions of Lemma \ref{lemma-for-def-temp-on-non-affine}. Denote the space of all tempered functions on $X$ by $\mathcal{T}(X)$.

\subsection{Proposition}\label{temp.-times-Schwartz-is-Schwartz_non_affine} Let $X$
be an algebraic variety, $t\in\mathcal{T}(X)$ and $s\in\mathcal{S}(X)$.
Then $t\cdot s\in\mathcal{S}(X)$.

Proof: Let $X=\bigcup\limits_{i=1}^{n}X_i$ be some affine cover, such that $s=\sum\limits_{i=1}^{n}Ext_{X_i}^{X}(s_i)$ for some $s_i\in\mathcal{S}(X_i)$. Then $t|_{X_i}\in\mathcal{T}(X_i)$ and by Proposition \ref{temp.-times-Schwartz-is-Schwartz} $t_i|_{X_i}\cdot s_i\in\mathcal{S}(X_i)$. Thus, $t\cdot s=\sum\limits_{i=1}^{n}Ext_{X_i}^{X}(s_i\cdot t|_{X_i})\in\mathcal{S}(X)$. \qed

\subsection{Definition}\label{def-X-flat-for-non-affine} Let $X$ be an algebraic variety. A function $f:X\to\bR$ is called \emph{flat at $x\in X$} if there exists an affine open neighborhood $x\in X_i \subset X$ such that $f|_{X_i}$ is flat at $x$. It is called flat at $Z\subset X$ if it is flat at any $x\in Z$.

Remark: Equivalently, a function $f:X\to\bR$ is called flat at $x\in X$ if {\bf for any} affine open neighborhood $x\in X_i \subset X$ one has that $f|_{X_i}$ is flat at $x$. This easily follows by intersecting any two affine open neighborhoods.

\subsection{Proposition (extension by zero for non affine varieties)}\label{ext-by-zero-for-non-affine}
Let $X$ be an algebraic variety, and $U$ an open subset of $X$. Then the extension by zero to $X$ of a Schwartz function on $U$ is a Schwartz function on $X$, which is flat at $X\setminus U$.

Proof: Consider some affine open cover $X=\bigcup\limits_{i=1}^kX_i$. Then $U=\bigcup\limits_{i=1}^k(U\cap X_i)$ is an affine open cover of $U$. Take some $s\in\mathcal{S}(U)$. By definition $s=\sum\limits_{i=1}^k Ext_{U\cap X_i}^{U}(s_i)$, for some $s_i\in\mathcal{S}(U\cap X_i)$. As $U\cap X_i$ is open in $X_i$, then, by Proposition \ref{ext-by-zero}, $Ext_{U\cap X_i}^{X_i}(s_i)$ is a Schwartz function on $X_i$, which is flat at $X_i\setminus (U\cap X_i)$. Then $Ext_{U}^{X}(s)=Ext_{U}^{X}(\sum\limits_{i=1}^k Ext_{U\cap X_i}^{U}(s_i))=\sum\limits_{i=1}^k Ext_{U}^{X}(Ext_{U\cap X_i}^{U}(s_i))=\sum\limits_{i=1}^k Ext_{X_i}^{X}(Ext_{U\cap X_i}^{X_i}(s_i))$ is by definition a Schwartz function on $X$, and clearly it is flat on $X\setminus U$. \qed

\subsection{Lemma (restrictions of tempered functions to closed and to open subsets for non affine varieties)}\label{temp.-rest-for-non-affine} Let $X$ be an algebraic variety, and let $U\subset X$ be some open subset. Then $Res_X^U(\mathcal{T}(X))\subset\mathcal{T}(U)$ and $Res_X^{X\setminus U}(\mathcal{T}(X))\subset\mathcal{T}(X\setminus U)$.

Proof: Consider some affine open cover $X=\bigcup\limits_{i=1}^kX_i$. Then $U=\bigcup\limits_{i=1}^k(U\cap X_i)$ is an affine open cover of $U$. Let $t\in\mathcal{T}(X)$, then by definition for any $1\leq i\leq k$, $t|_{X_i}\in\mathcal{T}(X_i)$. By Proposition \ref{Tempered-is-a-sheaf},  $(t|_{X_i})|_{U\cap X_i}\in\mathcal{T}(U\cap X_i)$. Clearly $(t|_{X_i})|_{U\cap X_i}=(t|_{U})|_{U\cap X_i}$, thus $t|_{U}\in\mathcal{T}(U)$, i.e. $Res_X^U(\mathcal{T}(X))\subset\mathcal{T}(U)$.

Observe that $X\setminus U=\bigcup\limits_{i=1}^k((X\setminus U)\cap X_i)$ is an affine open cover of $X\setminus U$. By Lemma \ref{temp.-rest}, $(t|_{X_i})|_{(X\setminus U)\cap X_i}\in\mathcal{T}((X\setminus U)\cap X_i)$. Clearly $(t|_{X_i})|_{(X\setminus U)\cap X_i}=(t|_{X\setminus U})|_{(X\setminus U)\cap X_i}$, and so $Res_X^{X\setminus U}(\mathcal{T}(X))\subset\mathcal{T}(X\setminus U)$.

\subsection{Corollary}\label{can-restrict-tempered-for-non-affine} Let $X$ be an algebraic variety, and let $U,V\subset X$ be two open subsets of $X$ such that $U\subset V$. Then $Res_V^U(\mathcal{T}(V))\subset \mathcal{T}(U)$.

Proof: $V$ is an algebraic variety, so we may apply Lemma \ref{temp.-rest-for-non-affine} to the algebraic variety $V$ and its open subset $U$. \qed

\subsection{Proposition}\label{Tempered-is-a-sheaf-for-non-affine} Let $X$ be an algebraic variety. The assignment of the space of tempered functions to any open $U\subset X$, together with the restriction of functions, form a sheaf on $X$.

Proof: By Corollary \ref{can-restrict-tempered-for-non-affine} the above is a pre-sheaf.

Let $\{U_i\}_{i\in I}$ be some open cover of $U$, and $t,s\in\mathcal{T}(U)$ be such that for any $i\in I$, $t|_{U_i}=s|_{U_i}$. As $t$ and $s$ are simply real valued functions on $U$, clearly $s=t$, and so the axiom of uniqueness holds.

Now let $t_i\in\mathcal{T}(U_i)$ be such that for any $i,j\in I$, $t_i|_{U_i\cap U_j}=t_j|_{U_i\cap U_j}$. Clearly there exists a (unique) function $t:U\to\bR$ such that for any $i\in I$, $t|_{U_i}=t_i$. In order to prove that the existence axiom holds, it is thus left to show that $t\in\mathcal{T}(U)$. By Proposition \ref{Zariski-is-Noetherian} we may assume $|I|<\infty$ by choosing some subcover and showing $t|_{U_i}=t_i$ only for indices $i$ in this subcover (as the functions we begin with agree on the intersections, this will automatically hold for all the other indices we omitted). By standard induction on the number of indices (i.e. the number of sets in the chosen finite subcover), it is enough to show that the following holds:

Let $X$ be an algebraic variety and let $U_1,U_2\subset X$ be two open subsets. Assume that for any $i\in\{1,2\}$ we are given $t_i\in \mathcal{T}(U_i)$ such that $t_1|_{U_1\cap U_2}=t_2|_{U_1\cap U_2}$. Then, there exists a function $t\in\mathcal{T}(U_1\cup U_2)$ such that $t|_{U_1}=t_1,t|_{U_2}=t_2$.

Clearly there exists a (unique) function $t:U_1\cup U_2\to\bR$ such that $t|_{U_i}=t_i$ for $1\leq i\leq 2$. It is left to show that $t\in\mathcal{T}(U_1\cup U_2)$. Consider some affine open cover $X=\bigcup\limits_{j=1}^kX_i$. Then $U_i=\bigcup\limits_{j=1}^k(U_i\cap X_j)$ is an affine open cover of $U_i$, and $U_1\cup U_2=\bigcup\limits_{j=1}^k((U_1\cup U_2)\cap X_j)$ is an affine open cover of $U_1\cup U_2$. As $t_i\in \mathcal{T}(U_i)$, one has $t_i|_{U_i\cap X_j}\in \mathcal{T}(U_i\cap X_j)$. As $(U_1\cup U_2)\cap X_j$ is affine, and $\bigcup\limits_{i=1}^2U_i\cap X_j$ is an affine open cover of it, and as $t_1|_{U_1\cap U_2}=t_2|_{U_1\cap U_2}$, then, by Proposition \ref{Tempered-is-a-sheaf}, $t|_{(U_1\cup U_2)\cap X_j}\in\mathcal{S}((U_1\cup U_2)\cap X_j)$, i.e. $t\in\mathcal{T}(U_1\cup U_2)$. \qed

{\bf Further work.} In order to prove that the rest of the properties that were proved in the affine case also hold in the general case, the next natural step should be proving a tempered partition of unity also holds in the non-affine case (i.e. to prove a non-affine version of Proposition \ref{temp-part-unity-affine-case}). Moreover, it seems that proving such a proposition would pave the way to prove all other properties. Our attempts of proving this in the algebraic context were not fruitful, however we suggest the following idea: as the theory of Schwartz spaces etc. is now fully established for arbitrary Nash manifolds and for affine real algebraic varieties (and partially established for non-affine real algebraic varieties), one may try to construct this theory for some bigger category, such that the Nash category and the algebraic category form subcategories of this category. The first guess, of taking the semi-algebraic category to be the nominated, will not work, as this category is in a sense "too flexible". As isomorphisms are not necessarily even smooth, in this category smooth functions are not pulled back by isomorphisms to smooth functions, and so there is no hope to define Schwartz spaces. We are currently working on generalizing the above theory to a category whose affine (local) models are closed semi-algebraic sets, but not all semi-algebraic morphisms are allowed. However, this category still generalizes both the Nash category and the algebraic category. If indeed such a theory can be established for this category, our hope (and guess) is that then proving partition of unity for non-affine varieties will be easier (as a consequence of the presence of "more" morphisms and "more" open subsets). Then, the results for the non-affine algebraic case would just follow as a special case.

\begin{appendix}
\section{Subanalytic geometry and proofs of two key lemmas }\label{subanalytic-geometry-2}\label{subanalytic-geometry-1}
The goal of this appendix is
to prove Lemmas \ref{X-flat-is-invariant}
and \ref{Milman_main_lemma}. All necessary preliminary needed results are
also given. In this appendix we always consider the Euclidean
topology on $\bR^n$, unless otherwise stated.

\subsection{Semi-analytic and subanalytic sets}\label{sem-sub-analy-set-def}

One may define semi-analytic sets in $\bR^n$ by (locally) using analytic
functions, namely: $A\subset\bR^n$ is semi-analytic if and only if for every
point $p\in A$ there exist an open neighbourhood $p\in U \subset \bR^n$ and
finitely many analytic functions $f_{i,j},g_{i,k}:U\to\mathbb{R}$, such that:
$$A\cap U=\bigcup_{i=1}^r\{x\in\bR^n|\forall 1\leq j \leq s_i,1\leq k\leq
t_i:f_{i,j}(x)>0,g_{i,k}(x)=0\}.$$ The images of semi-analytic sets under
analytic maps (which we did not define), and even under standard projections,
are not necessarily semi-analytic (see \cite[Example in the end of section
III]{L}).

This motivates the following definition (following \cite[Definition
3.1]{BM1}): $A\subset\bR^n$ is \emph{subanalytic} if and only if for every
point $p\in \bR^n$ there exists an open neighbourhood $p\in U \subset \bR^n$
such that $A\cap U$ is a projection of a relatively compact (i.e. bounded)
semi-analytic set (and see equivalent definitions in \cite[page 40 and page
95]{Shi2}). A map $\nu:A\to B$ (where $A\subset\bR^n,B\subset\bR^m$ are subanalytic)
is called \emph{subanalytic} if its graph is a subanalytic set in $\bR^{n+m}$.

\subsection{Definition}\label{flat_very_new_def} Let $X\subset\bR^n$ be a subanalytic set, $y\in X$ be some point and $m\in\mathbb{N}$. A function $f:X\to \bR$
is  \emph{$m$-flat at $y$} if there exists $F\in C^m(\bR^n)$ with $f=F|_X$,
such that the Taylor polynomial of order $m$ of $F$ at $y$ is identically
zero. If $f$ is $m$-flat at $y$ for any $m\in\mathbb{N}$, we will say that
$f$ is \emph{flat at $y$} (e.g. a function $f:\bR^n\to\bR$ is flat at some
$y\in\bR^n$ if the Taylor series of $f$ at $y$ is identically zero). If $f$
is ($m$-)flat at $y$ for any $y\in Z$ (where $Z\subset X$ is some subset),
we will say that $f$ is \emph{($m$-)flat at $Z$} \footnote{Defining $m$-flat
functions (and all other notions of flatness) on (closed) subanalytic sets
may be also done by two other different approaches. The first is by defining
the "Zariski $m$-paratangent bundle" and then defining $m$-flat functions
as restrictions of $C^m$ functions on the ambient space, that satisfy some
natural condition expressed in terms of the Zariski $m$-paratangent bundle.
The second is an intrinsic approach, using Glaeser extensions of real valued
functions. The interested reader is referred to \cite{BMP2,F}. It follows
from Theorems 1.7, 1.8 and 1.9 of \cite{BMP2} that all three approaches are
equivalent.}.

\subsubsection{An important remark} As the Taylor polynomial is only dependent
on the Euclidean local behaviour
of functions, one may substitute $\bR^n$ above by any Euclidean open subset
of $\bR^n$
containing $X$. This will be done in \ref{appendix-res-from-open-nbrhd},
and will be used when it will be more convenient.

\subsubsection{Equivalence of Definitions \ref{flat_new_def} and \ref{flat_very_new_def}} Any algebraic set is also subanalytic, and so it seems we have two different definitions of flat functions at a point: clearly if a function is flat at some point according to Definition \ref{flat_new_def} it is also flat according to Definition \ref{flat_very_new_def}. The other way is not trivial: by Definition \ref{flat_very_new_def} $f$ is flat
at $y$ means that for any $m\in\mathbb{N}$ there exists $F^m\in C^m(\bR^n)$
such that $F^m|_X=f$ and the Taylor polynomial of order $m$ of $F^m$ at $y$
is identically zero. This does not mean,
a-priori, that there exists $F\in C^{\infty}(\bR^n)$ such that $F|_X=f$ and
the Taylor series of $F$ at $y$ is identically zero, i.e. that the $f$ is flat at $y$ according to Definition \ref{flat_new_def}. According to \cite{M}, as $X$ is algebraic it is formally semicoherent rel.
the singleton $\{y\}$ (see \cite[Definition 1.2]{BM3} and discussion immediately after), which
is equivalent, according to \cite[Theorem 1.13]{BM3}, to the fact that $f$ does extend to such $F$ (later we will formulate it as $C^{\infty}(X;\{y\})=C^{(\infty)}(X;\{y\})$. Thus the two definitions of flat functions at a point in algebraic sets are equivalent.

\subsection{Restrictions from open neighborhoods and composite functions}\label{appendix-res-from-open-nbrhd}
In this subsection we will mainly follow the notations of \cite{BMP1}. Let
$X\subset\bR^n$ be some subanalytic set, $Z\subset X$ a closed subset (in
the induced topology from $\bR^n$), and $k\in\mathbb{N}$. As said, flatness
at a point is clearly a local property. Hence, it is natural to present the
following spaces of functions: $$C^k(X;Z):=\{f:X\to\bR|\exists U\subset\bR^n\text{
an open neighborhood of }X\text{ and }F\in C^k(U)$$$$\text{ such that }F|_X=f\text{
and }F\text{ is }k\text{-flat at }Z\}.$$

For $k=\infty$ we define similarly: $$C^\infty(X;Z):=\{f:X\to\bR|\exists
U\subset\bR^n\text{ an open neighborhood of }X\text{ and }F\in C^\infty(U)$$$$\text{
such that }F|_X=f\text{ and }F\text{ is flat at }Z\}.$$

Denote $C^{(\infty)}(X;Z):=\bigcap\limits_{k\in\mathbb{N}}C^k(X;Z)$. Clearly
$C^{\infty}(X;Z)\subset C^{(\infty)}(X;Z)$. In general $C^{\infty}(X;Z)\neq
C^{(\infty)}(X;Z)$, remarkably even when $Z=\emptyset$ (see, \cite{P}, Example
2).

For any $k\in\mathbb{N}\cup\{\infty\}$, denote $C^k(X):=C^k(X;\emptyset)$
-- this coincides with the usual definition if $X$ is smooth.

Let $\Omega$ be some open subset of $\bR^m$. A continuous map $\varphi:\Omega\to\bR^n$
is called semi-proper if for each compact subset $K$ of $\bR^n$ there exists
a compact subset $L$ of $\Omega$ such that $\varphi(L)=K\cap\varphi(\Omega)$.
Let $\varphi:\Omega\to \bR^n$ be a semi-proper real analytic map.  In that
case $X:=\varphi(\Omega)$ is a closed subanalytic subset of $\bR^n$. For
any $k\in\mathbb{N}$ define: $$(\varphi^*C^k(X))^\wedge:=\{f\in C^k(\Omega)|\forall
a\in X\text{ }\exists g\in C^k(X),f-\varphi^*(g):=f-g\circ \varphi\text{
is }k\text{-flat at }\varphi^{-1}(a)\},$$ and for $k=\infty$ define: $$(\varphi^*C^\infty(X))^\wedge:=\{f\in
C^\infty(\Omega)|\forall a\in X\text{ }\exists g\in C^\infty(X),f-\varphi^*(g):=f-g\circ
\varphi\text{ is flat at }\varphi^{-1}(a)\}.$$

Finally, for any closed subanalytic subset $Z\subset X$ and any $k\in\mathbb{N}\cup\{\infty\}$
define: $$(\varphi^*C^k(X;Z))^\wedge:=(\varphi^*C^k(X))^\wedge\cap C^k(\Omega;\varphi^{-1}(Z)).$$

We will use some properties of these spaces of functions, proved mainly in
\cite{BMP1}. Of special importance will be the following Theorem:

\subsubsection{Theorem (Uniformization theorem -- \cite[Theorem 5.1]{BM1})}\label{uniformization-theorem}
Let $M$ be a real analytic manifold, and let $X\subset M$ be a closed analytic
subset. Then there is a real analytic manifold $N$ and a proper real analytic
mapping $\tilde \varphi:N\to M$ such that $\tilde \varphi (N)=X$.

The following Lemma is stated in \cite{BMP1}, for the reader's convenience
we give here a detailed proof \footnote{The
authors thank Prof. Pierre D. Milman for providing guidelines for the proof
of Lemma \ref{Milman-easy-to-see-lemma}.}:

\subsubsection{Lemma}\label{Milman-easy-to-see-lemma} $\bigcap\limits_{k\in\mathbb{N}}(\varphi^*
C^k(X;Z))^{\wedge}=(\varphi^* C^\infty(X;Z))^{\wedge}$.

Proof: The inclusion $\supset$ is clear.

For the inclusion $\subset$ notice that $\bigcap\limits_{k\in\mathbb{N}}C^k(\Omega;\varphi^{-1}(Z))=C^\infty(\Omega;\varphi^{-1}(Z))$
as $\Omega$ is open (and in particular smooth). Thus, it is left to show
that $\bigcap\limits_{k\in\mathbb{N}}(\varphi^* C^k(X))^{\wedge}\subset(\varphi^*
C^\infty(X))^{\wedge}$. Fix some $f\in\bigcap\limits_{k\in\mathbb{N}}(\varphi^*
C^k(X))^{\wedge}$, $x\in X$, and $\omega\in\varphi^{-1}(x)$. Define the set
of $k$-jets at $x$ whose pullback$\text{'}$s $k$-jet at $\omega$ equals to
the $k$-jet of $f$, i.e. $$A^k:=\{J\text{ is a }k\text{-jet at }x|\varphi^*
J\text{'s }k\text{-jet at }\omega\text{ equals the }k\text{-jet of }f\text{
at }\omega\}.$$

By abuse of notation we considered jets as functions, i.e. chose a representative
-- this is independent of the choice made. Then the condition above may be
reformulated as $\varphi^*J-f$ is $k$-flat at $\omega$.

$A^k\neq\emptyset$, by the definition of $f$ (e.g. the $k$-jet of $g$ corresponding
to $x$ in the definition of $(\varphi^*C^k(X))^\wedge$ is contained in $A^k$).
For any $\mathbb{N}\cup\{0\}\ni l\leq k$ define the projection of $A^k$ to
$A^l$ by $A_l^k:=pr_l(A^k)\subset A^l$ (i.e. $A_l^k$ is the the set of all
$l$-jets that can be extended to $k$-jets in a way that is compatible with
$f$). Also define $A_l:=\bigcap\limits_{k\geq l}A_l^k\subset A^l$ (i.e. $A_l$
is the the set of all $l$-jets that can be extended to $k$-jets in a way
that is compatible with $f$ for any $k\geq l$). Assume the following hold
for any $l\in\mathbb{N}$:

(1) $A_l\neq\emptyset$.

(2) $pr_l:A_{l+1}\to A_l$ is onto, i.e. any $l$-jet in $A_l$ can be extended
to an $(l+1)$-jet in $A_{l+1}$.

In that case using the principle of dependent choices (see \ref{prelim-dependent-choice} below)
there exists a series of jets $\{J_l\}_{l\in\mathbb{N}}$ such that $J_l\in
A_l$ and for any $l\in\mathbb{N}$, $pr_l(J_{l+1})=J_l$. This series can be
thought of as a formal power series on $\bR^n$ (which we denote by $G$),
where $J_l$ is the truncated power series up to order $l$. By Borel's Theroem
(see \ref{Borel's-thm} below) there exists a function $\tilde G\in C^\infty(\bR^n)$
such that $T^\infty_x \tilde G=G$, where $T^\infty_x\tilde G$ means the Taylor
series of $\tilde G$ at $x$. By construction we observe that $T^\infty_\omega(\varphi^*\tilde
G)=T^\infty_\omega f$, and so conclude that $f\in(\varphi^*C^\infty(X))^\wedge$
(as $\tilde G|_X\in C^\infty(X)$ satisfies the desired condition).

Thus, it is left to prove that the assumptions (1) and (2) hold.

In order to see that $A_l\neq\emptyset$ observe that for any $k\geq l$, $A_l^k$
is a finite dimensional affine space (it is finite dimensional as it lies
in the finite dimensional space $A^l$). Moreover, these affine spaces form
a decreasing sequence $A^l=A_l^l\supset A_l^{l+1}\supset A_l^{l+2}\supset\cdot\cdot\cdot$.
As for any $k\geq l$, $A^k\neq \emptyset$, also $A_l^k\neq\emptyset$. The
dimensions of these affine spaces form a decreasing sequence of non-negative
integers, and so stabilizes. Thus we conclude that $\bigcap\limits_{k\geq
l}A_l^k(=A_l)=A_l^{k_l}\neq\emptyset$, for some $k_l\in\mathbb{N}$. And so
assumption (1) holds.

In order to see that $pr_l:A_{l+1}\to A_l$ is onto, we take some $J\in A_l$,
and define for any $k\geq l+1$: $B^k_{l+1}:=\{H\in A^k_{l+1}|pr_l(H)=J\}$.
By the definition of $A_l$ there exists a $k$-jet $J_k\in A^k$ such that
$pr_l(J_k)=J$, and so $pr_{l+1}(J_k)\in A_{l+1}^k$ satisfies $pr_l(pr_{l+1}(J_k))=J$,
and so $B^k_{l+1}\neq\emptyset$. As $B_{l+1}^{l+1}\supset B_{l+1}^{l+2}\supset
B_{l+1}^{l+3}\supset \cdot\cdot\cdot$ is a decreasing sequence of finite
dimensional affine spaces it stabilizes, and so $\bigcap\limits_{k\geq l+1}B_{l+1}^k\neq\emptyset$.
Choose some $\mathcal{J}\in \bigcap\limits_{k\geq l+1}B_{l+1}^k$. By definition
$\mathcal{J}\in A_{l+1}$ and by construction $pr_l(\mathcal{J})=J$. Thus,
assumption (2) holds. \qed

\subsubsection{The Principle of Dependent Choices (DC)}\label{prelim-dependent-choice}
We assume a weak version of the axiom of choice, the so called principle
of dependent choices: if $E$ is a binary relation on a nonempty set $A$,
and if for every $a\in
A$ there exists an element $b\in A$ such that $bEa$, then there is a sequence
$a_0,a_1,a_2,...,a_n,...$
in $A$, such that $a_{n+1}Ea_n$ for any $n\in\mathbb{N}$. The interested
reader is referred to \cite[page 50]{J}.

\subsubsection{Theorem (E. Borel -- \cite[Theorem 38.1]{T})}\label{Borel's-thm} Let $\Phi$
be an arbitrary formal power series in $n$ indeterminates, with complex coefficients.
Then there exists a function in $C^\infty (\bR^n)$ whose Taylor expansion
at the origin is identical to $\Phi$.

\subsection{Proof of Lemma \ref{Milman_main_lemma}}\label{new_label_for_milman} Recall we want to prove
the following: let $X$ be a compact (in the Euclidean topology) algebraic
set in $\bR^n$, and let $Z\subset X$ be some (Zariski) closed subset. Define
$U:=X\setminus Z$, $$W_Z:=\{\phi:X\to\bR|\exists\tilde\phi\in C^\infty(\bR^n)\text{
such that }\tilde\phi|_{X}=\phi\text{ and }\phi\text{
is }X\text{-flat on }Z\}$$ and $$(W^{\bR^n}_Z)^{comp}:=\{\phi\in C^\infty(\bR^n)|\phi\text{
is a compactly supported and is flat on }Z\}.$$ Then, for any $f\in W_Z$, there
exists $\tilde{f}\in (W_{Z}^{\bR^n})^{comp}$ such that $\tilde{f}|_{X}=f$.

We start the proof with a preliminary lemma:

\subsubsection{Lemma}\label{semi-proper-from-open}
There exists a pair $(\Omega,\varphi)$ such that $\Omega\subset\bR^m$ is
an open subset of $\bR^m$ (in the Euclidean topology), and $\varphi:\Omega\to\bR^n$
is a semiproper real analytic function, such that $\varphi(\Omega)=X$. Moreover,
for any $f\in C^{\infty}(X;\emptyset)$, $\varphi^*f:=f\circ \varphi\in C^{\infty}(\Omega)$.

Proof: By Theorem \ref{uniformization-theorem} there is a real analytic manifold
$N$ and a proper real analytic mapping $\tilde \varphi:N\to \bR^n$ such that
$\tilde \varphi (N)=X$. As $N$ is a real analytic manifold it has an open
cover $\{N_i\}_{i\in I}$ such that for any $i\in I$, $N_i$ is analytically
diffeomorphic to an open subset of $\bR^{d_i}$, i.e. there exist analytical
diffeomorphisms $\nu_i:N_i\to\bR^{d_i}$ such that $\nu_i(N_i)$ is open in
$\bR^{d_i}$. As $\tilde\varphi$ is proper and $X\subset\bR^n$ is compact,
$N=\tilde\varphi^{-1}(X)$ is compact, and thus there exists a finite subcover
$\{N_i\}_{i=1}^k$. We denote $m:=\max\limits_{i=1}^k\{d_i\}+1$. For any $1\leq
i\leq k$ define $\Omega_i=\nu_i(N_i)\times(i-\frac{1}{4},i+\frac{1}{4})^{m-d_i}\subset
\bR^m$, and define $\psi_i:\Omega_i\to \Omega_i$ by $\psi_i(n_i,\alpha_1,..,\alpha_{m-d_i}):=(n_i,i,..,i)$,
where $n_i\in \nu_i(N_i)$ and $\alpha_j\in (i-\frac{1}{4},i+\frac{1}{4})$.
Note that $\psi_i$ is semiproper real analytic. Define $\Omega:=\bigcup\limits_{i=1}^{k}\Omega_i$.
As $\Omega_1,..,\Omega_k$ are disjoint sets in $\bR^m$ we can naturally define
a semiproper function $\psi:\Omega\to\Omega$ by $\psi|_{\Omega_i}:=\psi_i$.
Clearly $\Omega\subset\bR^m$ is open. Now define a function $\nu^{-1}:\psi(\Omega)\to
N$ by $\nu^{-1}(n_i,i,..,i)=\nu_i^{-1}(n_i)$, where $n_i\in \nu_i(N_i)$.
Note that $\nu^{-1}$ is a proper map. Finally, defining $\varphi:=\tilde\varphi\circ\nu^{-1}\circ\psi$,
we get that $\varphi$ is a semiproper real analytic function, satisfying
$\varphi(\Omega)=X$. The Moreover part of Lemma \ref{semi-proper-from-open}
is obvious. \qed

Proof of Lemma \ref{Milman_main_lemma}:

Fix $f\in W_Z$. In particular , $f\in C^{\infty}(X;\emptyset)$. Then, by Lemma
\ref{semi-proper-from-open}, $\varphi^*f:=f\circ \varphi\in C^{\infty}(\Omega)$.
Denote $\tilde f:=\varphi^*f$. By definition $\tilde f\in (\varphi^*(C^{\infty}(X)))^{\wedge}$.

We will now prove that $\tilde f\in C^{\infty}(\Omega;\varphi^{-1}(Z))$. For any $z\in Z$ and any $k\in\mathbb{N}$ there exists $\bar f_z^k\in C^k(\bR^n)$
such that $f=\bar f_z^k|_X$ and $\bar f_z^k$ is $k$-flat at $z$. Note that
for any $z\in Z$ and any $k\in\mathbb{N}$, $\tilde f(=f\circ \varphi)=\bar
f_z^k\circ \varphi$, and in particular for any $z'\in Z$ and any $k'\in\mathbb{N}$
we have $\bar f_{z'}^{k'}\circ \varphi=\bar f_z^k\circ \varphi$. First we
will prove that for any $k\in\mathbb{N}$, $\tilde f\in C^{k}(\Omega;\varphi^{-1}(Z))$.
Fix some $z\in Z$ and some $\tilde z\in\varphi^{-1}(z)$.  As we may write
$\tilde f=\bar f_z^k\circ \varphi$, the fact that $\tilde f$ is $k$-flat
at $\tilde z$ follows immediately from Lemma \ref{Faa-Di-Bruno-for-Milman}.

We thus showed that for any $k\in\mathbb{N}$, $\tilde f\in C^{k}(\Omega;\varphi^{-1}(Z))$.
As we also had $\tilde f\in (\varphi^*(C^{k}(X)))^{\wedge}$, we get that
$\tilde f\in (\varphi^* C^k(X;Z))^{\wedge}:=(\varphi^*(C^{k}(X)))^{\wedge}\cap
C^{k}(\Omega;\varphi^{-1}(Z))$. As this holds for any $k\in\mathbb{N}$ we
get that $\tilde f\in \bigcap\limits_{k\in\mathbb{N}}(\varphi^* C^k(X;Z))^{\wedge}$.
By Lemma \ref{Milman-easy-to-see-lemma} $\bigcap\limits_{k\in\mathbb{N}}(\varphi^*
C^k(X;Z))^{\wedge}=(\varphi^* C^\infty(X;Z))^{\wedge}$, and by \cite[Theorem
1.3]{BMP1}, $(\varphi^* C^\infty(X;Z))^{\wedge}=\varphi^*C^{(\infty)}(X;Z)$,
so $\tilde f\in \varphi^*C^{(\infty)}(X;Z)$.

As $\tilde f\in \varphi^*C^{(\infty)}(X;Z)$ there exists $h\in C^{(\infty)}(X;Z)$
such that $\tilde f=\varphi^* h = h\circ \varphi$. Since $\tilde f =f\circ
\varphi$ and $\varphi$ is onto $X$ it follows that $h=f$, i.e. $f\in C^{(\infty)}(X;Z)$.
According to \cite{M}, as $X$ is algebraic it is formally semicoherent rel.
$Z$ (see \cite[Definition 1.2]{BM3} and discussion immediately after), which
is equivalent, according to \cite[Theorem 1.13]{BM3}, to the fact that $C^{\infty}(X;Z)=C^{(\infty)}(X;Z)$.
We conclude that $f\in C^{\infty}(X;Z)$.

Recall that we started with some $f\in W_Z$, we showed that $f\in C^{\infty}(X;Z)$,
and our goal is to show that $f$ is the restriction of some function in $(W_{Z}^{\bR^n})^{comp}$.
As $f\in C^{\infty}(X;Z)$ there exists an open $V\subset \bR^n$ such that
$X\subset V$, and $F:V\to\mathbb{R}$ such that $F\in C^\infty(V)$, $F$ is
flat at $Z$ and $F|_X=f$. Without loss of generality, as $X$ is compact in
the Euclidean topology on $\bR^n$, we may assume $V$ is bounded in the Euclidean
norm on $\bR^n$.

Take some open $V'\subsetneq V$ containing $X$. Let $\rho\in C^\infty(\bR^n)$
be a function supported in $V'$, such that $\rho|_{V''}=1$, where $V''\subsetneq
V'$ is some open subset containing $X$ (it is standard to show such $\rho$
exists by convolving the characteristic function of some open subset containing
$V''$ and strictly contained in $V'$, with some appropriate approximation
of unity). Now define $\tilde F:\bR^n\to\mathbb{R}$ by $\tilde F|_V:=\rho\cdot
F$ and $\tilde F|_{\bR^n\setminus V}:=0$. Clearly $\tilde F|_X=F|_X=f$, $\tilde
F\in C^\infty(\bR^n)$, $\tilde F$ is flat at $Z$ (as $\tilde F|_{V''}=F|_{V''}$),
and $\tilde F$ is compactly supported (in the Euclidean topology) in $\bR^n$,
i.e. $\tilde F\in (W_{Z}^{\bR^n})^{comp}$ and $\tilde F|_X=f$. \qed

\subsection{Multivariate Fa\'a Di Bruno Formula}\label{prelim-faa-di}  The
famous chain rule for deriving real valued functions from the real line stated
that $(f\circ g)'(x)=f'(g(x))\cdot g'(x)$. This can be generalized to higher
derivatives and higher dimensions, i.e. partial derivatives of arbitrary
order of composite multivariate functions. We will be interested only in
the following result:

\subsubsection{Lemma (cf. \cite[Theorem 2.1]{CS})}\label{Faa-Di-Bruno-for-Milman}
Let $x_0\in\bR^d$, $V\subset\bR^d$ be some open neighborhood of $x_0$ and
$g:V\to\bR^m,g\in C^k(V,\bR^m)$, for some $k\in\mathbb{N}$. Let $U\subset\bR^n$
be some open neighborhood of $g(x_0)$ and $f:U\to\mathbb{R},f\in C^k(U)$.
Assume $f$ is $k$-flat at $g(x_0)$, i.e. its Taylor polynomial of degree
$k$ at $g(x_0)$ is zero. Then $f\circ g:g^{-1}(U)\to\mathbb{R}$ is $k$-flat
at $x_0$.

\subsection{Proof of Lemma \ref{X-flat-is-invariant}}\label{proof-of-X-flat-is-inv.}

By definition for any $x\in X_1$, we have $\varphi(x)=(\frac{f_1(x)}{g_1(x)},..,\frac{f_{n_2}(x)}{g_{n_2}(x)})$,
where $f_1,..,f_{n_2},g_1,..,g_{n_2}\in\bR[X_1]$ and for any $1\leq i\leq
{n_2}$, $g_i^{-1}(0)\cap X_1=\emptyset$. By abuse of notation we choose some
representatives in $\bR[x_1,..,x_{n_1}]$ and consider $f_1,..,f_{n_2},g_1,..,g_{n_2}$
as functions in $\bR[x_1,..,x_{n_1}]$. Define $U:=\{x\in\bR^{n_1}|\prod\limits_{i=1}^{n_2}
g_i(x)\neq 0\}$. $U$ is open in $\bR^{n_1}$, $X_1$ is a closed subset of
$U$, and $\varphi$ can be naturally extended to a regular map $\tilde{\varphi}:U\to\bR^{n_2}$
(by the same formula of $\varphi$). Note that $U$ is an affine Nash manifold.

Let $f:X_2\to \bR$ be some function that is flat at some $p\in X_2$.
In particular $f\in C^{(\infty)}(X_2;\{p\})$. According to \cite{M}, as $X_2$
is algebraic it is formally semicoherent rel. to the singleton $\{p\}$ (see
\cite[Definition 1.2]{BM3} and discussion immediately after), which is equivalent,
according to \cite[Theorem 1.13]{BM3}, to the fact that $C^{\infty}(X_2;\{p\})=C^{(\infty)}(X_2;\{p\})$.
We conclude that $f\in C^{\infty}(X_2;\{p\})$, and so there exist an open
subset $V\supset X_2$ and a function $F\in C^{\infty}(V)$ such that $F$ is
flat at $p$ and $F|_{X_2}=f$. As $\tilde \varphi:U\to\bR^{n_2}$ is continuous,
$U':=\tilde \varphi^{-1}(V)$ is an open subset of $U$, and so an open subset
of $\bR^{n_1}$. Clearly $U'$ contains $X_1$. Denote by $G$ the pullback of
$F$ to $U'$ by $\tilde\varphi$, i.e. $G:U'\to\bR$ is defined by $G:=F\circ
(\tilde \varphi|_{U'})$. By Lemma \ref{Faa-Di-Bruno-for-Milman}, as
$F$ is flat at $p$, $G$ is flat at $\tilde\varphi^{-1}(p)$, and in particular
$f\circ \varphi=G|_{X_1}$ is flat at $\varphi^{-1}(p)$. \qed

\subsection{Remark}\label{remark-on-supp} Similar arguments to these given
in the first part of \ref{proof-of-X-flat-is-inv.} show that the notion of
support of a function on an affine algebraic variety $X$ is well defined.
Namely, we consider the support defined by {\bf some Euclidean topology}
on $X$, where we choose some closed embedding of $X$ in an affine space.
This support is independent of the embedding chosen.

\section{Cosheaves}\label{sheaves-and-cosheaves}

We recall the definition of a cosheaf on a topological space. For simplicity we assume our cosheves take values in the category of real vector spaces, but this can be replaced by any other Abelian category with arbitrary coproducts.

A \emph{pre-cosheaf} $F$ on a topological space $X$ is a covariant functor from $Top(X)$, where $Top(X)$ is the category whose objects are the open sets of $X$, and whose morphisms are the inclusion maps, to the category of real vector spaces. Equivalently, a pre-cosheaf on a topological space $X$ is the assignment of a real vector space $F(U)$ to any open $U\subset X$, and for any inclusion of open sets $V\subset U$ an extension morphism $Ext_V^U:F(V)\to F(U)$, such that $Ext_U^U=Id$ and for open $W\subset V\subset U$, $Ext_V^U\circ Ext_W^V=Ext_W^U$.

A \emph{cosheaf} on a topological space $X$ is a pre-sheaf, such that for any open $U\subset X$ and any open cover $\{U_i\}_{i\in I}$ of $U$, the following sequence is exact:

\begin{center}
\leavevmode
\xymatrix{
\bigoplus\limits_{i> j\in I}F(U_i\cap U_j) \ar[rr]^{Ext_1} & & \bigoplus\limits_{i\in I}F(U_i) \ar[rr]^{Ext_2} & & F(U) \ar[rr] & & 0,}
\end{center}

where $Ext_1(\bigoplus\limits_{i\neq j\in I}\xi_{i,j}):=\sum\limits_{i>j\in I}(Ext_{U_i\cap U_j}^{U_i}(\xi_{i,j})-Ext_{U_i\cap U_j}^{U_j}(\xi_{i,j}))$, and $Ext_2(\bigoplus\limits_{i\in I}\xi_i):=\sum\limits_{i\in I}Ext^{U}_{U_i}(\xi_i)$. Equivalently, a cosheaf on a topological space $X$ is a pre-sheaf, such that for any open $U\subset X$ and any open cover $\{U_i\}_{i\in I}$ of $U$, the following two axioms hold:

\begin{enumerate}
\item For any $s\in F(U)$ there exists a finite set of indexes $i_1,..,i_k$ such that for any $1\leq j\leq k$ there exists $s_{i_j}\in F(U_{i_j})$, and $\sum\limits_{j=1}^k Ext_{U_{i_j}}^U(s_{i_j})=s$ (this corresponds to fact that $Ext_2$ is onto).
\item For any finite subset $J\subset I$, if we are given for any $i\in J$ some element $s_i\in F(U_i)$, such that $\sum\limits_{i\in J}Ext_{U_i}^{U}s_i=0$, then there exists a finite subset $J\subset J'\subset I$, such that for any $i> j\in J'$ there exists $s_{i,j}\in F(U_i\cap U_j)$ satisfying $s_i=\sum\limits_{J' \ni j<i}Ext_{U_i\cap U_j}^{U_i}(s_{i,j})-\sum\limits_{J' \ni j>i}Ext_{U_i\cap U_j}^{U_i}(s_{j,i})$ for any $i\in J$ (this corresponds to the fact that $Im(Ext_1)=Ker(Ext_2)$).
\end{enumerate}

A cosheaf on a topological space $X$ is said to be \emph{flabby} if for any two open subsets $U,V\subset X$ such that $V\subset U$, the extension morphism $Ext_V^U:F(V)\to F(U)$ is injective.

\section{Noetherianity of the Zariski topology}\label{proof-of-Notherianity}
First, as by definition any algebraic variety has a {\bf finite} cover by
affine varieties, it is enough to prove Proposition \ref{Zariski-is-Noetherian}
for affine varieties. Second, it is enough to prove Proposition \ref{Zariski-is-Noetherian}
for algebraic subsets of $\bR^n$. Let $X\subset \bR^n$ be an algebraic set,
$U\subset X$ a Zariski open subset, and $\{U_\alpha\}_{\alpha\in I}$ an open
cover of $U$. We prove Proposition \ref{Zariski-is-Noetherian} in 3 steps:

Case 1 -- Assume $X=U=\bR^n$. By Proposition \ref{any-alg-set-is-a-zero-locus}
for any $\alpha\in I$ there exists $f_\alpha\in\bR[x_1,..,x_n]$ such that
$U_\alpha=\{x\in\bR^n|f_{\alpha}(x)\neq 0\}$. As $\bR^n\setminus \bigcup\limits_{\alpha\in
I}U_\alpha=\emptyset$, the zero locus of $<f_\alpha>_{\alpha\in I}$ (the
set of all points $x\in\bR^n$ satisfying $f(x)=0$ for any $f$ in the ideal
generated by all of the polynomials $\{f_\alpha\}_{\alpha\in I}$) is empty.
By Hilbert's basis Theorem $\bR[x_1,..,x_n]$ is Noetherian, and so there
exist $g_1,..,g_m\in \bR[x_1,..,x_n]$ such that $<f_\alpha>_{\alpha\in I}=<g_1,..,g_m>$.
As the zero locus of this ideal is empty, $g:=\sum\limits_{l=1}^m g_l^2$
satisfies $g^{-1}(0)=\emptyset$. As $g\in <f_\alpha>_{\alpha\in I}$ there
exist $a_1,..,a_k\in \bR[x_1,..,x_n]$ and $\alpha_1,..,\alpha_k\in I$ such
that $g=\sum\limits_{i=1}^k a_i\cdot f_{\alpha_k}$. This implies that $f_{\alpha_1},..,f_{\alpha_k}$
have no common zeroes, and so $\bigcup\limits_{i=1}^k U_{\alpha_i}=\bR^n$.

Case 2 -- Assume $X=U$. There exist $\{V_\alpha\}_{\alpha\in I}$ open
in $\bR^n$ such that $U_\alpha=V_\alpha\cap X$. Then $\{V_\alpha\}_{\alpha\in
I}\cup (\bR^n\setminus X)$ is an open cover of $\bR^n$. By case (1) it has
a finite subcover, and so intersecting this subcover with $X$ we get a finite
subcover of $U$.

Case 3 -- The general case. By proposition \ref{open-of-alg-is-aff},
$U$ is itself an affine algebraic variety. Moreover, $\{U_\alpha\cap U\}_{\alpha\in
I}$ is an open cover of $U$. By considering some closed embedding of $U$
in some $\bR^m$, we are reduced to case (2). \qed

\end{appendix}

\

\end{document}